\documentclass[12pt]{amsart}
\usepackage{amsmath,amssymb,setspace,nicefrac, yhmath, amscd,eucal,enumitem}
\usepackage[active]{srcltx}

\allowdisplaybreaks
\usepackage[all]{xy}

\newtheorem{theorem}{Theorem}
\newtheorem{lemma}{Lemma}%
\newtheorem{corollary}{Corollary}%

\newtheorem{problem}{Problem}%

\theoremstyle{thmstyletwo}%
\newtheorem{example}{Example}%
\newtheorem{remark}{Remark}%

\theoremstyle{thmstylethree}%
\newtheorem{definition}{Definition}%

\newcommand{\cE}{{\mathcal E}}

\newcommand{\cM}{{\mathcal M}}
\newcommand{\cP}{{\mathcal P}}

\newcommand{\cR}{{\mathcal R}}

\newcommand{\cH}{{\mathcal H}}

\newcommand{\cF}{{\mathcal F}}

\newcommand{\cD}{{\mathcal D}}

\newcommand{\N}{{\mathbb N}}
\newcommand{\minus}{\backslash}

\newcommand{\Adm}{{\rm Adm}}

\newcommand{\HK}{{\rm HK}}
\newcommand{\ET}{{\rm ET}}
\newcommand{\WLET}{{\rm WLET}}

\begin{document}

\title[Weak Optimal Entropy Transport Problems]{Weak Optimal Entropy Transport Problems}


\author{Nhan-Phu Chung}
\address{Nhan-Phu Chung, Institute of Applied Mathematics, University of Economics Ho Chi Minh City, Vietnam.}
\email{phuchung82@gmail.com; phucn@ueh.edu.vn} 

\author{Thanh-Son Trinh}
\address{Thanh-Son Trinh, Faculty of Information Technology, Industrial University of Ho Chi Minh City, Vietnam.}
\email{trinhthanhson@iuh.edu.vn} 







\begin{abstract}In this paper, we introduce weak optimal entropy transport problems that cover both optimal entropy transport problems and weak optimal transport problems introduced by Liero, Mielke, and Savar\'{e} \cite{LMS18}; and Gozlan, Roberto, Samson and Tetali \cite{GRST}, respectively. Under some mild assumptions of entropy functionals, we establish a Kantorovich type duality for our weak optimal entropy transport problem. As consequences, via a different method, we recover both Kantorovich duality formulas for optimal entropy transport problems \cite{LMS18}, and weak optimal transport problems \cite{GRST,BBP}.
\end{abstract}




\maketitle
\section{Introduction}
After pioneering works of Kantorovich in 1940s \cite{Kant42,Kant48}, the theory of classical Monge-Kantorovich optimal transport problems has been developed by many authors. It has many applications in other fields such as economics, geometry of nonsmooth metric spaces, image processing, PDEs, functional inequalities, probability and statistics,... We refer to the monographs \cite{AGS,Galichon,PC,San,V03,V} for a more detailed presentation and references therein. The primal Monge-Kantorovich problem is written in the form 
\begin{align*}
\inf\bigg\{\int_{X_1\times X_2} cd\boldsymbol{\gamma}:\boldsymbol{\gamma}\in \Pi(\mu_1,\mu_2)\bigg\},
\end{align*}
where $\mu_1, \mu_2$ are given probability measures on Polish metric spaces $X_1$ and $X_2$, $c:X_1\times X_2\to (-\infty,+\infty]$ is a cost function, and $\Pi(\mu_1,\mu_2)$ is the set of all probability measures $\boldsymbol{\gamma}$ on $X_1\times X_2$ with marginals $\mu_1$ and $\mu_2$. 
	
Recently, in a seminal paper \cite{LMS18}, Liero, Mielke and Savar\'{e} introduced theory of Optimal Entropy Transport problems between nonnegative and finite Borel measures in Polish spaces which may have different masses. Since then it has been investigated further in \cite{CP,CT,D,FMS, KV,LM,MS}.
They relaxed the marginal constraints $\gamma_i:=\pi^i_\sharp \boldsymbol{\gamma}=\mu_i$ via adding penalizing divergences
\begin{align*}
\cF_i(\gamma_i\vert\mu_i):=\int_{X_i}F_i(f_i(x_i))d\mu_i(x_i)+(F_i)'_\infty \gamma_i^\perp(X),
\end{align*}   
where $\gamma_i=f_i\mu+\gamma_i^\perp$ is the Lebesgue decomposition of $\gamma_i$ with respect to $\mu_i$, and $F_i:[0,\infty)\to [0,\infty]$ are given convex, lower semi-continuous functions with their recession constants $(F_i)'_\infty:=\lim_{s\to \infty}\dfrac{F_i(s)}{s}$. Such functions will be referred to as entropy functions in the sequel. Then the Optimal Entropy Transport problem is formulated as 
\begin{align}
\label{P-OET}
\mathcal{ET}(\mu_1,\mu_2):=\inf_{\boldsymbol{\gamma}\in \cM(X_1\times X_2)} \cE(\boldsymbol{\gamma}\vert \mu_1,\mu_2),
\end{align}
where $\cE(\boldsymbol{\gamma}\vert\mu_1,\mu_2):=\sum_{i=1}^2\cF_i(\gamma_i\vert \mu_i)+\int_{X_1\times X_2}c(x_1,x_2)d\boldsymbol{\gamma}(x_1,x_2),$ and $\cM(X_1\times X_2)$ is the space of all nonnegative and finite Borel measures on $X_1\times X_2$.
Given entropy functions $F_1, F_2:[0,\infty)\to [0,\infty]$, we define functions $F_i^\circ:\mathbb{R}\to [-\infty,\infty]$ and $R_i:[0,\infty)\to [0,\infty]$ by $F_i^\circ(\varphi):=\inf_{s\geq 0}(\varphi s+F_i(s))$ for every $\varphi\in \mathbb{R}$, and
\begin{align}
	R_i(r):=\left\{\begin{array}{ll}
		rF(1/r) &\text{ if }r>0,\\
		(F_i)'_{\infty} &\text{ if }r=0.
	\end{array}\right.\notag
\end{align}

In \cite{LMS18}, the authors showed that under certain mild conditions of entropy functions $F_i$, the problem \eqref{P-OET} always has minimizing solutions and they established the following duality formula
\begin{align*}
\mathcal{ET}(\mu_1,\mu_2)&=\sup_{(\varphi_1,\varphi_2)\in \boldsymbol{\Phi}}\sum_{i=1}^2\int_{X_i}F_i^\circ(\varphi_i)d\mu_i\\ 
&=\sup_{(\psi_1,\psi_2)\in \boldsymbol{\Psi}}\sum_{i=1}^2\int_{X_i}\psi_id\mu_i,
\end{align*}
where 
\begin{align*}
\boldsymbol{\Phi}&:=\bigg\{(\varphi_1,\varphi_2)\in C_b(X_1,\mathring{D}(F_1^\circ))\times C_b(X_2,\mathring{D}(F_2^\circ)):\varphi_1\oplus \varphi_2\leq c\bigg\}, \\
\boldsymbol{\Psi}&:=\bigg\{(\psi_1,\psi_2)\in C_b(X_1,\mathring{D}(R_1^*))\times C_b(X_2,\mathring{D}(R_2^*)):R_1^*(\psi_1)\oplus R_2^*(\psi_2)\leq c\bigg\}.
\end{align*}
Here $\mathring{D}(F)$ is the interior of $D(F):=\{r\geq 0: F(r)<\infty\}$, $f_1\oplus f_2\leq c$ means that $f(x_1)+f_2(x_2)\leq c(x_1,x_2)$ for every $x_1\in X_1, x_2\in X_2$,	$C_b(A,B)$ is the set of all continuous and bounded functions from $A$ to $B$, $F^*:  \mathbb{R} \to (-\infty,+\infty]$ is the Legendre conjugate function of $F$ defined by 
\begin{align*}
F^*(\varphi):=\sup_{s\geq 0}(s\varphi-F(s)) \mbox{ for every } \varphi\in \mathbb{R}.
\end{align*}	

On the other hand, in 2014, Gozlan, Roberto, Samson and Tetali \cite{GRST} introduced weak optimal transport problems encompassing the classical Monge-Kantorovich optimal transport and weak transport costs introduced by Talagrand and Marton in the 90's. After that, theory of weak optimal transport problems and its applications have been investigated further by numerous authors \cite{ABC,ACJ,BBHK, BBP,BP,GT,GJ,GRSST, Shu}. In \cite{GRST}, the authors also established a Kantorovich type duality for their weak optimal transport problem as follows. 

Let $\cP(X_2)$ be the space of all Borel probability measures on $X_2$ and $C:X_1\times \cP(X_2)\to [0,\infty]$ be a lower semi-continuous function such that $C(x,\cdot)$ is convex for every $x\in X_1$. Given $\mu_1\in \cP(X_1),\mu_2\in \cP(X_2)$ and $\boldsymbol{\gamma}\in \Pi(\mu_1,\mu_2)$, we denote its disintegration with respect to the first marginal $\gamma_1$ by $(\gamma_{x_1})_{x_1\in X_1}$. Then  
	the weak optimal transport problem is defined as
	\begin{align}
	\label{P-weak optimal transport}
	V(\mu_1,\mu_2):=\inf\bigg\{\int_{X_1}C(x_1,\gamma_{x_1})d\mu_1(x_1):\boldsymbol{\gamma}\in \Pi(\mu_1,\mu_2)\bigg\}, 
	\end{align}
	and its Kantorovich duality is
	$$ V(\mu_1,\mu_2)=\sup\bigg\{\int_{X_1}R_C\varphi(x_1)d\mu_1(x_1)-\int_{X_2}\varphi(x_2)d\mu_2(x_2):\varphi\in C_b(X_2)\bigg\},$$
	where $$R_C\varphi(x_1):=\inf_{p\in \cP(X_2)}\bigg\{\int_{X_2}\varphi(x_2)dp(x_2)+C(x_1,p)\bigg\}, \text{ for all }x_1\in X_1.$$

In this paper, we introduce weak optimal entropy transport (WOET) problems which generalize both optimal entropy transport \cite{LMS18} and weak optimal transport problems \cite{GRST}. For every $\boldsymbol{\gamma}\in \cM(X_1\times X_2)$, we denote by $\gamma_1,\gamma_2$ the first and second marginals of $\boldsymbol{\gamma}$. We also denote its disintegration with respect to the first marginal $\gamma_1$ by $(\gamma_{x_1})_{x_1\in X_1}$ i.e, for every bounded Borel function $f: X_1\times X_2\to   \mathbb{R} $ we have  
	$$\int_{X_1\times X_2}fd\boldsymbol{\gamma} =\int_{X_1}\left(\int_{X_2}f(x_1,x_2)d\gamma_{x_1}(x_2)\right)d\gamma_1(x_1),$$
	where $\gamma_1$ is the first marginal of $\boldsymbol{\gamma}$. 
	Given $\mu_1\in \cM(X_1), \mu_2\in \cM(X_2)$, our primal weak optimal entropy transport problem is formulated as
\begin{align}
\label{P-primal}
\cE_C(\mu_1,\mu_2):=\inf_{\boldsymbol{\gamma}\in \cM(X_1\times X_2)} \cE_C(\boldsymbol{\gamma}\vert\mu_1,\mu_2) , 
\end{align}
		where 
\begin{align}	
\label{D-relative entropy functions}	
		\cE_C(\boldsymbol{\gamma}\vert\mu_1,\mu_2):=\sum_{i=1}^2\cF_i(\gamma_i\vert\mu_i)+\int_{X_1}C(x_1,\gamma_{x_1})d\gamma_1(x_1).
\end{align}		
Before stating the main results of the article, let us introduce some notations. Let $F_i:[0,\infty)\to [0,\infty]$, $i=1,2$ be admissible entropy functions. 
We define	
\begin{align*}
	\Lambda:=\bigg\{(\varphi_1,\varphi_2)\in C_b(X_1,\mathring{ D}(F_1^\circ))\times C_b(X_2,\mathring{ D}(F_2^\circ))&:\varphi_1(x_1)+p(\varphi_2)\leq C(x_1,p),\\
 &\mbox{ for every }
	 x_1\in X_1, p\in \cP(X_2)\bigg\},
\end{align*}
and 
\begin{align*}
\Lambda_R:=\bigg\{&(\varphi_1,\varphi_2)\in C_b(X_1)\times C_b(X_2):\sup_{x_i\in X_i}\varphi_i(x_i)<F_i(0),i=1,2,\\
& \mbox{ and } R_1^*(\varphi_1(x_1))+p\left(R^*_2(\varphi_2)\right)\leq C(x_1,p) \mbox{ for every } x_1\in X_1,p\in \cP(X_2)\bigg\}. 
\end{align*}
Our main result is a Kantorovich duality for our weak optimal entropy transport problem.	
\begin{theorem}\label{dual theorem}
		Let $X_1,X_2$ be Polish metric spaces. Let $C:X_1\times \cP(X_2)\to (-\infty,+\infty]$ be a lower semi-continuous function such that $C$ is bounded from below and $C(x_1,\cdot)$ is convex for every $x_1\in X_1$. Let $F_i:[0,\infty)\to [0,\infty]$, $i=1,2$ be admissible entropy functions such that $F_i$ is superlinear, i.e. $(F_i)'_\infty=+\infty$ for $i=1,2$. Then for every $\mu_i\in \cM(X_i),i=1,2$ we have that 
\begin{align*}
\cE_C(\mu_1,\mu_2)&=\sup_{(\varphi_1,\varphi_2)\in \Lambda}\sum_{i=1}^2\int_{X_i}F_i^\circ(\varphi_i)d\mu_i\\
&=\sup_{(\varphi_1,\varphi_2)\in \Lambda_R}\sum_{i=1}^2\int_{X_i}\varphi_id\mu_i\\
&=\sup_{\varphi\in C_b(X_2)}\int F_1^\circ (R_C \varphi)d\mu_1+\int F_2^\circ (-\varphi) d\mu_2.
\end{align*}	 	
	\end{theorem}
Assume that there exists some cost function $c:X_1\times X_2\to (-\infty,+\infty]$ which is lower semi-continuous and bounded from below, such that $C(x_1,p)=\int_{X_2}c(x_1,x_2)dp(x_2)$ for every $x_1\in X_1,p\in \cP(X_2)$ then our (WOET) problem (\ref{P-primal}) becomes the Optimal-Entropy Transport problem (\ref{P-OET}). Furthermore, in this case it is not difficult to check that $\Phi=\Lambda$ (Lemma \ref{L-temp1}) and $C$ is lower semi-continuous (Lemma \ref{L-lower semicontinuous in Martingale problems}). Therefore, via a different proof, we recover the duality formula of Optimal Entropy Transport problem in \cite[Theorem 4.11 and Corollary 4.12]{LMS18} when $F_1, F_2$ are superlinear.
\begin{corollary}
Let $X_1,X_2$ be Polish metric spaces. Let $c:X_1\times X_2\to (-\infty,+\infty]$ be a lower semi-continuous function which is bounded from below. Let $F_i:[0,\infty)\to [0,\infty]$, $i=1,2$ be admissible entropy functions such that $F_i$ is superlinear, i.e. $(F_i)'_\infty=+\infty$ for $i=1,2$. Then for every $\mu_i\in \cM(X_i),i=1,2$ we have that 
\begin{align*}
\mathcal{ET}(\mu_1,\mu_2)&=\sup_{(\varphi_1,\varphi_2)\in \boldsymbol{\Phi}}\sum_{i=1}^2\int_{X_i}F_i^\circ(\varphi_i)d\mu_i.
\end{align*}
\end{corollary}

On the other hand, if we consider the admissible entropy functions $F_1,F_2:[0,\infty)\rightarrow [0,\infty]$ defined by $$F_1(r)=F_2(r):=\left\lbrace\begin{array}{ll}
		0 & \text{ if } r=1,\\
		+\infty &\text{ otherwise,}
	\end{array}\right.$$	
	then given $\mu_1\in \cM(X_1),\mu_2\in \cM(X_2)$, our (WOET) problem  will become the following pure weak transport problem. 
	\begin{equation}
\cE_C(\mu_1,\mu_2)=\inf_{\boldsymbol{\gamma}\in \cM(X_1\times X_2)}\left\{\int_{X_1}C(x_1,\gamma_{x_1})d\gamma_1(x_1) \vert \pi^i_\sharp\gamma=\mu_i,i=1,2\right\}.\label{P-temp}
	\end{equation}
		In this example, if $\boldsymbol{\gamma}\in \cM(X_1\times X_2)$ is a feasible plan, i.e. there exists $\boldsymbol{\gamma}\in \cM(X_1\times X_2)$ such that $\cE_C(\boldsymbol{\gamma}\vert\mu_1,\mu_2)<\infty$ then $\mu_1,\mu_2$ are the marginals of $\boldsymbol{\gamma}$. Thus, a necessary condition for feasibility is that $\vert \mu_1\vert=\vert\mu_2\vert$. If furthermore $\mu_i\in \cP(X_i), i=1,2$ then \eqref{P-primal} will be the weak transport problem \eqref{P-weak optimal transport}. In this case, we have that $F_1, F_2$ are superlinear and $F_i^\circ(\varphi)=\inf_{s\geq 0}(s\varphi+F_i(s))=\varphi$ for every $\varphi\in \mathbb{R}$. 
Therefore, from Theorem \ref{dual theorem} we get the following corollary which recovers a Kantorovich duality formula of the weak optimal transport problem established in \cite{BBP,GRST}. 
\begin{corollary}
Let $X_1,X_2$ be Polish metric spaces. Let $C:X_1\times \cP(X_2)\to (-\infty,+\infty]$ be a lower semi-continuous function such that $C$ is bounded from below and $C(x_1,\cdot)$ is convex for every $x_1\in X_1$. Let $F_1,F_2:[0,\infty)\rightarrow [0,\infty]$ be admissible entropy functions defined by $$F_1(r)=F_2(r):=\left\lbrace\begin{array}{ll}
		0 & \text{ if } r=1,\\
		+\infty &\text{ otherwise,}
	\end{array}\right.$$
		 Then for every $\mu_1\in \cP(X_1),\mu_2\in \cP(X_2)$ we have that
		 \begin{align*}
		 \cE_C(\mu_1,\mu_2)&=V(\mu_1,\mu_2)\\
		 &=\sup_{\varphi\in C_b(X_2)}\int F_1^\circ (R_C \varphi)d\mu_1+\int F_2^\circ (-\varphi) d\mu_2\\
		 &=\sup\bigg\{\int_{X_1}R_C\varphi(x_1)d\mu_1(x_1)-\int_{X_2}\varphi(x_2)d\mu_2(x_2):\varphi\in C_b(X_2)\bigg\}.
		 		 \end{align*}  
\end{corollary} 	
On the other hand, for the case $X_1$ and $X_2$ are compact, using a different approach which is inspired from the proof of \cite[Theorem 4.11]{LMS18} we can relax superlinear condition of $F_1, F_2$ for our duality formula. However, we need to add an extra assumption that the primal problem is feasible.	
\begin{theorem}\label{assume compact} Assume that $ X_1,X_2 $ are compact and $ (F_1)'_\infty+(F_2)'_\infty+\inf C>0 $. Let $ \mu_1\in\cM(X_1),\mu_2\in\cM(X_2) $. If problem \eqref{P-primal} is feasible, i.e. there exists $\boldsymbol{\gamma}\in \cM(X_1\times X_2)$ such that 
$\cE_C(\boldsymbol{\gamma}\vert\mu_1,\mu_2)<\infty$ then we have $$\cE_C(\mu_1,\mu_2)=\sup_{(\varphi_1,\varphi_2)\in \Lambda}\sum_{i=1}^2\int_{X_i}F_i^\circ(\varphi_i)d\mu_i.$$  
	\end{theorem}	
\begin{remark}
To prove a Kantorovich duality in the classical optimal transport problems for general Polish metric spaces, we often prove for the compact case first. Then using it for compact subsets of the spaces and combining this with approximation processes we will get the result, see for example \cite[Section 1.3]{V03}. To establish a Kantorovich duality for the optimal entropy transport problems in Polish metric spaces in \cite[Theorem 4.11]{LMS18}, the authors also did in this way. However, as optimal entropy transport problems have penalizing divergences $\cF_1,\cF_2$, induced from entropy functions $F_1, F_2$, this process is more complicated than the classical case. For our (WOET) problems, we not only deal with penalizing divergences $\cF_1,\cF_2$ but also the disintegrations of marginals. The latter term makes this approximation process from compact cases to general cases challenging. Therefore, to prove Theorem \ref{dual theorem} we really need a different method from \cite{LMS18}.   
\end{remark}	

Let us describe our strategy to prove Theorem \ref{dual theorem}. The inequality 
$$\cE_C(\mu_1,\mu_2)\geq\sup_{(\varphi_1,\varphi_2)\in \Lambda}\sum_{i=1}^2\int_{X_i}F_i^\circ(\varphi_i)d\mu_i$$
is easy to establish, and we only need a mild condition that $ (F_1)'_\infty+(F_2)'_\infty+\inf C>0 $ to get it (Lemma \ref{L-easy part of the dual formula}). The difficult part is to prove the converse inequality
\begin{align}
\label{temp-difficulty inequality}
\cE_C(\mu_1,\mu_2)\leq\sup_{(\varphi_1,\varphi_2)\in \Lambda}\sum_{i=1}^2\int_{X_i}F_i^\circ(\varphi_i)d\mu_i.
\end{align}

Given a metric space $X$, we denote by $(C_b(X))^*$  the dual space of the normed space $(C_b(X),\|\cdot\|_\infty)$. For every $\mu\in \cM(X)$, the map $T_\mu:C_b(X)\to \mathbb{R}$, defined by $f\mapsto \int_X fd\mu$, is a bounded linear operator, i.e. it belongs to $(C_b(X))^*$. Now to prove \eqref{temp-difficulty inequality} we define the functional $\ET:(C_b(X_1))^*\times(C_b(X_2))^*\to [-\infty,+\infty]$ as follows
 \begin{align}
 \label{D-functional ET}
\ET(T_1,T_2):=\left\lbrace\begin{array}{ll}
		\cE_C(\mu_1,\mu_2) & \text{ if } (T_1,T_2)=(T_{\mu_1}, T_{\mu_2}) 
		,\\
		+\infty &\text{ otherwise.}
	\end{array}\right.
\end{align}
 Given $\mu,\nu\in \cM(X)$, if $\int_X fd\mu=\int_X fd\nu$ for every $f\in C_b(X)$ then one gets $\mu=\nu$ \cite[Theorem 5.9, page 39]{Par}. Therefore, for every metric space $X$ we can consider $\cM(X)$ as a subset of $(C_b(X))^*$, and hence the functional $\ET$ is well defined.

 For the convenience, we will write $\ET(\mu_1,\mu_2)$ for $\ET(T_{\mu_1}, T_{\mu_2})$	for every $(\mu_1,\mu_2)\in \cM(X_1)\times \cM(X_2)$. We define
\begin{align*}
	\Lambda_{ET} :=\bigg\{&(\varphi_1,\varphi_2)\in C_b(X_1)\times C_b(X_2):\sum_{i=1}^2\int_{X_i}\varphi_id\mu_i\leq \ET(\mu_1,\mu_2),\\
	& \mbox{ for every } (\mu_1,\mu_2)\in \cM(X_1)\times\cM(X_2)\bigg\}, 
	\end{align*}
	$$\Lambda_{ET}^<:=\{(\varphi_1,\varphi_2)\in \Lambda_{ET} \vert \sup_{x_i\in X_i}\varphi_i(x_i)<F_i(0),i=1,2\}.$$	
Then we show that \begin{align*}
\cE_C(\mu_1,\mu_2)&=\sup_{(\varphi_1,\varphi_2)\in \Lambda_{ET}}\sum_{i=1}^2\int_{X_i}F_i^\circ(\varphi_i)d\mu_i\\
&=\sup_{(\varphi_1,\varphi_2)\in \Lambda_{ET}^<}\sum_{i=1}^2\int_{X_i}F_i^\circ(\varphi_i)d\mu_i.
\end{align*}
After that, we prove $\Lambda_{ET}^<=\Lambda_R$ and the inequality (\ref{temp-difficulty inequality}). Our proof of Theorem 1 relies on the fact that the functional $\ET$, defined in \eqref{D-functional ET}, is convex and positively homogenous and lower semi-continuous, and is thus the support function of a convex set. This fact is established in Lemma \ref{L-lower semicontinuous of ET} and Lemma \ref{L-Lamda_{ET}^<=Lamda_R}. The same strategy has been used in the proof of Theorem 4.2 of the paper \cite{ABC} by Alibert-Bouchitt\'{e}-Champion, dealing with duality for Weak Optimal Transport problems.
   	 
Our paper is organized as follows. In section 2, we review notations and properties of entropy functionals. In section 3, we prove Theorem \ref{dual theorem} and Theorem \ref{assume compact}. In this section we also investigate the existence of minimizers and the feasibility of our (WOET) problems. 
Finally, we will illustrate examples of our results including the ones that cover optimal entropy transport problems \cite{LMS18}, weak optimal transport problems \cite{ABC, BBP, GRST}.   

In a companion paper \cite{CT1}, we study a weak optimal entropy transport problem in which the entropy functions $F_i$, $i=1,2$ are not superlinear. 

\textbf{Acknowledgement.} Part of this work has been done when N-P. Chung and T-S. Trinh were affiliated with Sungkyunkwan University, Korea, and we were supported by the National Research Foundation of Korea (NRF) grants funded by the Korea government No. NRF- 2016R1A5A1008055 and No. NRF-2019R1C1C1007107. N-P. Chung is also funded by University of Economics Ho Chi Minh City, Vietnam. We thank the two anonymous referees for useful comments, particular for suggesting both the statement of Lemma \ref{L-BM} and its proof, which helped us completely remove a technical condition in the previous manuscript and significantly improve our results.

	\section{Preliminaries}	
	Let $(X,d) $ be a metric space. We denote by $\mathcal{M}(X)$ (resp. $\cP(X)$) the set of all positive Borel measures (resp. probability Borel measures) with finite mass. We denote by $C_b(X)$ the space of all real valued continuous bounded functions on $X$. 
	
	For any $\mu\in \cM(X)$, set $\vert \mu\vert:=\mu(X)$. Let $M$ be a subset of $\cM(X)$. We say that $M$ is \textit{bounded} if there exists $C>0$ such that $\vert \mu\vert \leq C$ for every $\mu\in M$, and $M$ is \textit{equally tight} if for every $\varepsilon>0$, there exists a compact subset $K_\varepsilon$ of $X$ such that $\mu\left(X\backslash K_\varepsilon\right)\leq \varepsilon$ for every $\mu\in M$. 
	
A metric space $X$ is \textit{Polish} if it is complete and separable. The weak topology on $\cM(X)$ is the smallest topology such that for each $f\in C_b(X)$, the map $\mu\mapsto \int_Xfd\mu$ is continuous, i.e. a sequence $\{\mu_n\}_{n\in \N}\subset \cM(X)$ converges weakly to $\mu\in \cM(X)$ if and only if $\lim_{n\to \infty}\int_Xfd\mu_n=\int_Xfd\mu$ for every $f\in C_b(X)$. We recall Prokhorov's Theorem.
	\begin{theorem} (Prokhorov's Theorem)
		Let $(X,d)$ be a Polish metric space. Then a subset $M\subset \cM(X)$ is bounded and \textit{equally tight} if and only if $M$ is relatively compact under the weak topology. 
	\end{theorem}
	
	Let $\mu_1,\mu_2\in \cM(X)$. If $\mu_2(A)=0$ yields $\mu_1(A)=0$ for any Borel subset $A$ of $X$ then we say that $\mu_1$ is \textit{absolutely continuous} with respect to $\mu_2$ and write $\mu_1\ll \mu_2$. We call that $\mu_1\perp \mu_2$ if there exists a Borel subset $A$ of $X$ such that $\mu_1(A)=\mu_2(X\minus A)=0$.
	
	 Let $\mu,\gamma\in \cM(X)$ then there are a unique measure $\gamma^\perp\in \cM(X)$ and a unique $\sigma\in L^1_+(X,\mu)$ such that $\gamma=\sigma\mu+\gamma^\perp, \mbox{ and } \gamma^\perp \perp\mu$. It is called the \textit{Lebesgue decomposition} of $\gamma$ relative to $\mu$.  
	
	Let $X_1,X_2$ be metric spaces. For any $\boldsymbol{\gamma}\in\cM(X_1\times X_2)$, we call that $\gamma_1$ and $\gamma_2$ are the first and second marginals of $\boldsymbol{\gamma}$ if $$\boldsymbol{\gamma}(A_1\times X_2)=\gamma_1(A_1)\text{ and }\boldsymbol{\gamma}(X_1\times A_2)=\gamma_2(A_2),$$
	for every  Borel subsets $A_i$ of $X_i$, $i=1,2$. Given $\mu_1\in \cM(X_1),\mu_2\in \cM(X_2)$, we denote by $\Pi(\mu_1,\mu_2)$ the set of all Borel measures on $X_1\times X_2$ with marginals $\mu_1$ and $\mu_2$. It is clear that $\Pi(\mu_1,\mu_2)$ is nonempty if and only if $\mu_1$ and $\mu_2$ have the same masses.
	
	Let $f:X_1\to X_2$ be a Borel map and $\mu\in\cM(X_1)$. We denote by $f_\sharp \mu \in\cM(X_2)$ the \textit{push-forward measure} defined by $$f_\sharp \mu(B):=\mu(f^{-1}(B)),$$
	for every Borel subset $B$ of $X_2$.

	We now review on entropy functionals. For more details, readers can see \cite[Section 2]{LMS18}.	
	
	We define the class of \textit{admissible entropy functions} by \begin{align*}
		\Adm(  \mathbb{R} _+):=\{F:[0,\infty)\to[0,\infty]& \vert F\text{ is convex, lower semi-continuous }\\
		&\text{and }D(F)\cap (0,\infty)\neq \emptyset \},
	\end{align*}
	where $D(F):=\{s\in [0,\infty) \vert F(s)<\infty\}$. We also denote by $\mathring{D}(F)$ the interior of $D(F)$.
	
	Let $F\in \Adm(  \mathbb{R} _+)$, we define function $F^\circ:  \mathbb{R} \to [-\infty,\infty]$ by 
	\begin{align}
	\label{D-temp 5}
	F^\circ(\varphi):=\inf_{s\geq 0}\big(\varphi s+F(s)\big) \text{ for every } \varphi\in   \mathbb{R} .
	\end{align}
	
	Given $F\in \Adm(  \mathbb{R} _+)$ we define the recession constant $F'_\infty$ by
	\begin{align}
	\label{F-recession constant}
	F_\infty':=\lim_{s\to\infty}\frac{F(s)}{s}, 
	\end{align}
		and we define the functional $\cF:\cM(X)\times \cM(X)\to [0,\infty]$   by 
	$$\cF(\gamma \vert \mu):=\int_XF(f)d\mu+F'_\infty \gamma^\perp(X),$$
	where $\gamma=f\mu+\gamma^\perp$ is the Lebesgue decomposition of $\gamma$ with respect to $\mu$, and we adopt the convention that $\infty \cdot 0=0$. 	
	
	The Legendre conjugate function $F^*:  \mathbb{R} \to (-\infty,+\infty]$ is defined by 
\begin{align}
\label{D-temp 7}
F^*(\varphi):=\sup_{s\geq 0}(s\varphi-F(s)).
\end{align}	
	
	Then it is clear that $F^\circ(\varphi)=-F^*(-\varphi), $ for every $\varphi\in   \mathbb{R} $. Note that $\mathring{D}(F^*)=(-\infty,F'_\infty)$ and $F^*$ is continuous and non-decreasing on $(-\infty,F'_\infty)$ \cite[page 989]{LMS18} and hence we get that
	\begin{align}
	\label{temp 10}
	\mathring{D}(F^\circ)=(-F'_\infty,+\infty) \mbox{ and } F^\circ \mbox{ is non-decreasing on } (-F'_\infty,+\infty).
	\end{align}

	Next, we define the reverse density function $R:[0,\infty)\to [0,\infty]$ of a given $F\in \Adm(  \mathbb{R} _+)$ by 
\begin{align}
\label{D-temp 6}
	R(r):=\left\{\begin{array}{ll}
		rF(1/r) &\text{ if }r>0,\\
		F'_{\infty} &\text{ if }r=0.
	\end{array}\right.
\end{align}	
	It is not difficult to check that the function $R$  is convex, lower semi-continuous, and $R(0)=F'_\infty$, $R'_\infty=F(0)$. Then $R\in \Adm(  \mathbb{R} _+)$. From \cite[the first line, page 992]{LMS18} we have  
	\begin{align}
	\label{F-interior of R*}
	\mathring{D}(R^*)=(-\infty, F(0)).
	\end{align}

	We also define the functional $\cR:\cM(X)\times\cM(X)\to [0,\infty]$ by $$\cR(\mu \vert \gamma):=\int_XR(\varrho)d\gamma+R'_\infty\mu^\perp(X),$$
where $\mu=\varrho\gamma+\mu^\perp$ is the Lebesgue decomposition of $\mu$ with respect to $\gamma$.	

Then by \cite[Lemma 2.11]{LMS18} for every $\mu,\gamma\in \cM(X)$ we have that 
\begin{align}
\label{F-relation between entropy and its reverse}
\cF(\gamma \vert \mu)=\cR(\mu \vert \gamma).
\end{align}

\begin{lemma}
	(\cite[Lemma 2.6 and formula (2.17)]{LMS18})
	\label{L-equality of F for (BM)}
	Let $X$ be a Polish space, $\gamma,\mu\in \cM(X)$. Let $F\in \Adm(\mathbb{R}_+)$ and $\phi,\psi: X\to [-\infty,+\infty]$ be Borel functions such that 
\begin{enumerate}
	\item $\cF(\gamma\vert \mu)<\infty;$
	\item $\psi(x)\leq F^*(\phi(x)) \text{ if } -\infty<\phi(x)\leq F'_\infty,\phi(x)<+\infty,$
	\item	$\psi(x)=-\infty \text{ if } \phi(x)=F'_\infty=+\infty,$
	\item	$\psi(x)\in [-\infty,F(0)] \text{ if }  \phi(x)=-\infty.$
\end{enumerate}
			If $\psi_-\in L^1(X,\mu)$ (resp. $\phi_-\in L^1(X,\gamma)$) then $\phi_+\in L^1(X,\gamma)$ (resp. $\psi_+\in L^1(X,\mu)$) and 
	\begin{align}
	\label{F-Temp1}
	\cF(\gamma\vert\mu)-\int_X\psi d\mu\geq \int_X\phi d\gamma.
	\end{align}
	Assume further that $\psi\in L^1(X,\mu)$ or $\phi\in L^1(X,\gamma)$, and $\mu=\rho\gamma$ for some $\rho\in L^1(X,\gamma)$ with $\rho(x)>0$ for every $x\in X$. Then equality holds in \eqref{F-Temp1} if and only if $\phi(x)=-R^*(\psi(x))$, and 
	$$\rho(x)\in D(R),\mbox{ } \psi(x)\in D(R^*), \mbox{ } R(\rho(x))+R^*(\psi(x))=\rho(x)\psi(x),$$ for $\mu$-a.e in $X$.  
	\end{lemma}
	
\begin{lemma}\label{L-inequality of F}
	(\cite[Theorem 2.7 and Remark 2.8]{LMS18})
	Let $X$ be a Polish space, $\gamma,\mu\in \cM(X)$ and $F\in \Adm(  \mathbb{R} _+)$. Then
	\begin{align*}
	\cF(\gamma \vert \mu)&=\sup\bigg\{\int_X F^\circ (\varphi)d\mu-\int_X\varphi d\gamma:\varphi\in C_b(X,\mathring{ D}(F^\circ))\bigg\}\\
	&=\sup\bigg\{\int_X \psi d\mu-\int_X R^*(\psi)d\gamma:\psi\in C_b(X,\mathring{ D}(R^*))\bigg\}.
 \end{align*}	 
	\end{lemma}

	\section{Weak optimal entropy transport problems}
	Let $X_1,X_2$ be Polish spaces. For every $\boldsymbol{\gamma}\in \cM(X_1\times X_2)$, we denote its disintegration with respect to the first marginal by $(\gamma_{x_1})_{x_1\in X_1}$ i.e, for every bounded Borel function $f: X_1\times X_2\to   \mathbb{R} $ we have  
	$$\int_{X_1\times X_2}fd\boldsymbol{\gamma} =\int_{X_1}\left(\int_{X_2}f(x_1,x_2)d\gamma_{x_1}(x_2)\right)d\gamma_1(x_1),$$
	where $\gamma_1$ is the first marginal of $\boldsymbol{\gamma}$. Note that $\gamma_{x_1}$ is a Borel probability measure on $X_2$ for every $x_1\in X_1.$

		We consider a function $C:X_1\times \cP(X_2)\to (-\infty,+\infty]$ which is lower semi-continuous, bounded from below and satisfies that for every $x\in X_1, \mbox{ }C(x,\cdot)$ is convex, i.e.
\begin{align}
\label{Convexity of C}	
	C(x,tp+(1-t)q)\leq tC(x,p)+(1-t)C(x,q),
	\end{align}
	for every $t\in [0,1], p,q\in \cP(X_2)$. 
	
	Given $F_1,F_2\in \Adm(  \mathbb{R} _+)$ and $\mu_1\in \cM(X_1),\mu_2\in \cM(X_2)$, we investigate the following problem.
	\begin{problem}
	\label{P-Weak Optimal Entropy-Transport Problem}
	(Weak Optimal Entropy-Transport Problem) Find $\bar{\boldsymbol{\gamma}}\in \cM(X_1\times X_2)$ minimizing
		$$\cE_C(\bar{\boldsymbol{\gamma}} \vert \mu_1,\mu_2)=\cE_C(\mu_1,\mu_2):=\inf_{\boldsymbol{\gamma}\in \cM(X_1\times X_2)} \cE_C(\boldsymbol{\gamma} \vert \mu_1,\mu_2) \mbox{ (WOET)}, $$
		where $\cE_C(\boldsymbol{\gamma} \vert \mu_1,\mu_2):=\sum_{i=1}^2\cF_i(\gamma_i \vert \mu_i)+\int_{X_1}C(x_1,\gamma_{x_1})d\gamma_1(x_1),$ and $\gamma_1,\gamma_2$ are the first and second marginals of $\boldsymbol{\gamma}.$
	\end{problem}
	\begin{remark}
		As we will see in Examples \ref{E-Optimal Entropy-Transport Problems} and \ref{E-Weak Optimal Transport Problems} in section 4, our (WOET) problems cover the Optimal Entropy-Transport problem \eqref{P-OET} and the Weak Optimal Transport problem \eqref{P-weak optimal transport}.  	
	\end{remark}
	First, we investigate the feasibility of Problem \ref{P-Weak Optimal Entropy-Transport Problem}. We say that Problem \ref{P-Weak Optimal Entropy-Transport Problem} is feasible if there exists $\boldsymbol{\gamma}\in \cM(X_1\times X_2)$ such that $\cE_C(\boldsymbol{\gamma} \vert \mu_1,\mu_2)<\infty$. 
	\begin{lemma} 
	\label{L-feasible}
	Let $\mu_1\in \cM(X_1)$ and $\mu_2\in \cM(X_2)$ with $m_i:=\mu_i(X_i)$. Then 
	\begin{enumerate}
	\item If Problem \ref{P-Weak Optimal Entropy-Transport Problem} is feasible then $K\ne \varnothing $, where 
	$$ K:=\bigg(m_1D(F_1)\bigg )\cap\bigg(m_2D(F_2)\bigg);$$
	\item Problem \ref{P-Weak Optimal Entropy-Transport Problem} is feasible if one of the following conditions is satisfied
	\begin{enumerate}[label=(\roman*)]
	\item both $F_i(0)<\infty, i=1,2$;
	\item the set $K\neq \varnothing$, $m_1m_2\ne 0$, and there exist $ B_i\in L^1(X_i,\mu_i) $ for $i=1,2$ with \[ C(x_1,p)\le B_1(x_1)+p(B_2) \mbox{ for every } x_1\in X_1, p\in \cP(X_2).\]
	\end{enumerate}	
	\end{enumerate}
	\end{lemma}
\begin{proof}
(1) Let $\boldsymbol{\gamma}\in \cM(X_1\times X_2)$ such that $\cE_C(\boldsymbol{\gamma} \vert \mu_1,\mu_2)<\infty$. From \cite[(2.44)]{LMS18}, we have $ \cF_i(\gamma_i \vert \mu_i)\geq m_iF_i( \vert \gamma_{i} \vert /m_i) $ for every $ i $. Thus, $ m_iF_i( \vert \gamma_{i} \vert /m_i)<\infty $ for every $ i=1,2 $. Hence, $  \vert \boldsymbol{\gamma} \vert \in m_iD(F_i) $ for every $ i $ and therefore the set $ K $ is not empty.

(2) (i) Let $\boldsymbol{\gamma}_0\in \cM(X_1\times X_2)$ be the null measure. Then 
$$\cE_C(\mu_1,\mu_2)\leq  \cE_C(\boldsymbol{\gamma}_0 \vert \mu_1,\mu_2)\leq\sum_{i=1}^{2}F_i(0) \vert \mu_i \vert <\infty. $$ 
(ii) Considering the Borel measure $ \boldsymbol{\gamma}=\dfrac{\theta}{m_1m_2}\mu_1\otimes\mu_2 $ with $ \theta\in K.$ Then we have 
\begin{align*} \cE_C(\boldsymbol{\gamma} \vert \mu_1,\mu_2)=& m_1F_1(\theta/m_1)+m_2F_2(\theta/m_2)+\int_{X_1} C\left(x_1,\frac{1}{m_2}\mu_2\right)d\dfrac{\theta}{m_1}\mu_1\\
\le&m_1F_1(\theta/m_1)+m_2F_2(\theta/m_2)+\int_{X_1} B_1(x_1)+\frac{1}{m_2}\mu_2(B_2)d\dfrac{\theta}{m_1}\mu_1\\
\le&m_1F_1(\theta/m_1)+m_2F_2(\theta/m_2)+\sum\theta m_i^{-1}\ \vert B_i\ \vert _{L^1(X_i,\mu_i)}<\infty.
 \end{align*}
   \end{proof}	
	
Next, we will show the existence of minimizers of (WOET) problems under some mild assumptions on $F_i$. 	
	\begin{lemma}\label{L-lower semicontinuous of integration C}
		Let $\{\boldsymbol{\pi}^k\}\subset \cM(X_1\times X_2)$ such that $\boldsymbol{\pi}^k$ converges to $\boldsymbol{\pi}$ in the weak topology. Then $$\liminf_{k\rightarrow \infty}\int_{X_1}C(x_1,\pi_{x_1}^k)d\pi^k_1(x_1)\geq \int_{X_1}C(x_1,\pi_{x_1})d\pi_1(x_1).$$
	\end{lemma} 
	\begin{proof}
		If $ \{\boldsymbol{\pi}^k\}\subset\cP(X_1\times X_2) $ then the result follows from \cite[Proposition 2.8]{BBP}. Now we consider the general case $ \{\boldsymbol{\pi}^k\}\subset\cM(X_1\times X_2) $. Since $C$ is bounded from below, there exists $K\in  \mathbb{R} $ such that $\overline{C}(x_1,p):=C(x_1,p)+K\geq 0$ for every $x_1\in X_1$ and $p\in\cP(X_2)$. If $ \boldsymbol{\pi} $ is the null measure then \begin{align*}
\liminf_{k\rightarrow \infty}\int_{X_1}C(x_1,\pi_{x_1}^k)d\pi^k_1(x_1)=\liminf_{k\rightarrow \infty}\bigg(\int_{X_1}\overline{C}(x_1,p)d\pi^k_1(x_1)-K \vert \boldsymbol{\pi}^k \vert \bigg)\geq 0.		
\end{align*}
So that we get the inequality. Note that by weak convergence, $  \vert \boldsymbol{\pi}^k \vert =\int 1d\boldsymbol{\pi}^k\to\int 1d\boldsymbol{\pi}= \vert \boldsymbol{\pi} \vert  $. If $\boldsymbol{\pi} $ is not the null measure then for sufficient large index $ k $ we also have $ \boldsymbol{\pi}^k $ is not the null measure. For convenience, just consider $ \boldsymbol{\pi}^k $ is not the null measure for all $ k $. For any $ \varphi\in C_b(X_1\times X_2) $ we have \begin{align*} \left \vert \int\varphi\left(\dfrac{1}{ \vert \boldsymbol{\pi}^k \vert }-\dfrac{1}{ \vert \boldsymbol{\pi} \vert }\right)d\boldsymbol{\pi}^k\right \vert \leq\left \vert \dfrac{1}{ \vert \boldsymbol{\pi}^k \vert }-\dfrac{1}{ \vert \boldsymbol{\pi} \vert }\right \vert \ \vert \varphi\ \vert _\infty \vert \boldsymbol{\pi}^k \vert \to 0, 
		\end{align*}
		and $ \int\varphi\frac{1}{ \vert \boldsymbol{\pi} \vert }d\boldsymbol{\pi}^k\to\int\varphi\frac{1}{ \vert \boldsymbol{\pi} \vert }d\boldsymbol{\pi} $. Therefore, $ \frac{\boldsymbol{\pi}^k}{ \vert \boldsymbol{\pi}^k \vert } $ weakly converges to $ \frac{\boldsymbol{\pi}}{ \vert \boldsymbol{\pi} \vert } $. Applying the result of the case `$ \{\boldsymbol{\pi}^k\}\subset\cP(X_1\times X_2) $' we get \begin{align*}
			\liminf_{k\to\infty}\int_{X_1}C(x_1,\pi_{x_1}^k)d\pi^k_1(x_1)=&\lim_{k\to \infty} \vert \boldsymbol{\pi}^k \vert \liminf_{k\to\infty}\int_{X_1}C(x_1,\pi_{x_1}^k)d\frac{\pi^k_1}{ \vert \boldsymbol{\pi}^k \vert }(x_1)\\
			\geq& \vert \boldsymbol{\pi} \vert \int_{X_1}C(x_1,\pi_{x_1})d\frac{\pi_1}{ \vert \boldsymbol{\pi} \vert }(x_1)\\
			=&\int_{X_1}C(x_1,\pi_{x_1})d\pi_1(x_1).
		\end{align*}
	   \end{proof}	
	
	\begin{theorem}\label{T-existence of minimizers}
		Let $ \mu_1\in\cM(X_1),\mu_2\in\cM(X_2) $ such that the problem (WOET) is feasible. We also assume that one of the following conditions (coercive conditions) hold:
		\begin{enumerate}[label=\roman*)]
			\item the entropy functions $ F_1 $ and $ F_2 $ are superlinear, i.e. $ (F_1)'_\infty=(F_2)'_\infty=+\infty;$
			\item the spaces $ X_1 $ and $ X_2 $ are compact and $ (F_1)'_\infty+(F_2)'_\infty+\inf C>0. $
		\end{enumerate}
		Then, the problem (WOET) admits a minimizer.
	\end{theorem}
	\begin{proof} By Lemma \ref{L-lower semicontinuous of integration C} and \cite[Corollary 2.9]{LMS18}, we get that for every $\mu_i\in \cM(X_i),i=1,2$ the functional $\cE(\cdot \vert \mu_1,\mu_2)$ is lower semi-continuous in $\cM(X_1\times X_2)$. Let $ \boldsymbol{\gamma}^n\subset \cM(X_1\times X_2)$ be a minimizing sequence of the problem (WOET).  
		
		For the case i), as $ \cE(\boldsymbol{\gamma}^n \vert \mu_1,\mu_2) $ is bounded above, $ \cF_i $ is non-negative  and $ C $ is bounded from below we get that $ \cF_i(\gamma_i^n \vert \mu_i) $ is bounded above. Applying \cite[Proposition 2.10]{LMS18}, the set $ \{\boldsymbol{\gamma}^n_i\} $ is a subset of a bounded and equally tight set. Hence, so is $ \{\gamma^n_i\} $ for each $ i $ and so is $ \{\boldsymbol{\gamma}^n\} $ by \cite[Lemma 5.2.2]{AGS}.
		
		For the case ii), if one of $ (F_i)'_\infty>0 $ then by applying \cite[Proposition 2.10]{LMS18} $ \boldsymbol{\gamma}^n $ is bounded as $ \boldsymbol{\gamma}^n(X_1\times X_2)=\boldsymbol{\gamma}^n_i(X_i). $ We only need to consider $ (F_i)'_\infty=0 $ for every $ i. $ In that case, we have $ \boldsymbol{\gamma}^n(X_1\times X_2)\le\dfrac{1}{\inf C}\cE_C(\boldsymbol{\gamma}^n \vert \mu_1,\mu_2). $ So $ \{\boldsymbol{\gamma}^n \} $ is bounded.
		
		In both cases, $ \{\boldsymbol{\gamma}^n \} $ is relatively compact by Prokhorov's Theorem and the proof is complete.
	   \end{proof}
Now we will prove our duality formulations of the (WOET) problems.
	
	We recall  
\begin{align*}
	\Lambda:=\bigg\{(\varphi_1,\varphi_2)\in C_b(X_1,\mathring{ D}(F_1^\circ))\times C_b(X_2,\mathring{ D}(F_2^\circ))&:\varphi_1(x_1)+p(\varphi_2)\leq C(x_1,p),\\
 &\mbox{ for every }
	 x_1\in X_1, p\in \cP(X_2)\bigg\},
\end{align*}	
and 
\begin{align*}
\Lambda_R:=\bigg\{&(\varphi_1,\varphi_2)\in C_b(X_1)\times C_b(X_2):\sup_{x_i\in X_i}\varphi_i(x_i)<F_i(0),i=1,2,\\
& \mbox{ and } R_1^*(\varphi_1(x_1))+p\left(R^*_2(\varphi_2)\right)\leq C(x_1,p) \mbox{ for every } x_1\in X_1,p\in \cP(X_2)\bigg\}. 
\end{align*}
\begin{lemma}\label{L-the set lamda is nonempty}
Let $X_1,X_2$ be Polish spaces and assume that $(F_1)'_\infty+(F_2)'_\infty+\inf C>0$ then $\Lambda$ is a nonempty set. If moreover $F_i$ is superlinear for $i=1,2$ then $\Lambda_R$ is also nonempty.
\end{lemma}
\begin{proof}
We consider first the case one of $F_1,F_2$ is superlinear. Assume that $(F_2)'_\infty=+\infty$ then by \eqref{temp 10} we get $\mathring{D}(F^\circ_2)=  \mathbb{R} $. Since $C$ is bounded from below one get that there exists $K\in   \mathbb{R} $ such that $C(x_1,p)\geq K$ for every $x_1\in X_1,p\in\cP(X_2)$. Let $\varepsilon>0$ and putting $\varphi_1(x_1):=\varepsilon$ on $X_1$ and $\varphi_2(x_2):=K-\varepsilon$ on $X_2$. Then $\varphi_2\in C_b(X_2,\mathring{D}(F^\circ_2))$. From \eqref{temp 10} we get $\mathring{D}(F_1^\circ)=(-(F_1)'_\infty,+\infty)$, and  hence as $(F_1)'_\infty\geq 0$ we have $\varphi_1\in C_b(X_1,\mathring{D}(F^\circ_1))$. Furthermore, for every $x_1\in X_1,p\in\cP(X_2)$ one has
$$\varphi_1(x_1)+p(\varphi_2)=\varepsilon+K-\varepsilon=K\leq C(x_1,p).$$
Thus, $(\varphi_1,\varphi_2)\in\Lambda$.

If $(F_1)'_\infty=+\infty$ then using the same argument as above we still have $\Lambda$ is nonempty.

Next, we consider the case $(F_i)'_\infty<\infty$ for $i=1,2$. As $(F_1)'_\infty+(F_2)'_\infty+\inf C>0$, there is $a>0$ such that $\inf C> (-(F_1)'_\infty+a)+(-(F_2)'_\infty+a)$. Then set $\varphi_i(x_i):=-(F_i)'_\infty+a$ on $X_i$ for $i=1,2$. From this we have $\varphi_i\in C_b(X_i,\mathring{D}(F^\circ_i))$ for $i=1,2$ and $$\varphi_1(x_1)+p(\varphi_2)=(-(F_1)'_\infty+a)+(-(F_2)'_\infty+a)<\inf C\leq C(x_1,p),$$
for every $x_1\in X_1,p\in\cP(X_2)$. Therefore, $\Lambda$ is nonempty.

Now, we assume $F_i$ is superlinear, we will prove that $\Lambda_R$ is nonempty. Suppose that $C$ is bounded below by $2S$. As $\mathring{D}(F^\circ_i)=  \mathbb{R} $ one has $F^\circ_i(S)\in   \mathbb{R} $ for $i=1,2$. Fixed $\varepsilon>0$ and set $\varphi_i(x_i):=\min\{F^\circ_i(S),F_i(0)\}-\varepsilon$ on $X_i$ for $i=1,2$. Then it is clear that $\varphi_i\in C_b(X_i)$ and $\sup_{x_i\in X_i}\varphi_i(x_i)<F_i(0)$. Notice that we also get $\varphi_i(x_i)<F_i^\circ(S)=-F_i^*(-S)$ on $X_i$ for $i=1,2$. By \cite[(2.31)]{LMS18} we obtain that $R^*_i(\varphi_i(x_i))\leq S$ on $X_i$ for $i=1,2$. Hence, for every $x_1\in X_1,p\in\cP(X_2)$ one has
$$R^*_i(\varphi_1(x_1))+p(R^*_2(\varphi_2))\leq 2S\leq C(x_1,p).$$
This means that $(\varphi_1,\varphi_2)\in\Lambda_R$.

   \end{proof}
Now we prove Theorem \ref{assume compact}.
 \begin{proof}[Proof of Theorem \ref{assume compact}.]
		Put $M:=\{\boldsymbol{\gamma}\in\cM(X_1\times X_2) \vert 
		\int_{X_1}C(x_1,\gamma_{x_1})d\gamma_1(x_1)<\infty\}$ and $B:=C_b(X_1,\mathring{ D}(F_1^\circ))\times C_b(X_2,\mathring{ D}(F_2^\circ))$. As our primal problem (\ref{P-primal}) is feasible, we must have that $M$ is not an empty set. Let $\boldsymbol{\gamma},\boldsymbol{\bar{\gamma}}\in M$ and $t\in [0,1]$. By the convexity of $C(x_1,\cdot)$ and observe that $\bigg((1-t)(d\gamma_1/d((1-t)\gamma_1+t\bar{\gamma}_1)\gamma_{x_1}+t(d\bar{\gamma}_1/d((1-t)\gamma_1+t\bar{\gamma}_1)\bar{\gamma}_{x_1}\bigg)_{x_1\in X_1}$ is the disintegration of the measure $\boldsymbol{\beta}:=(1-t)\boldsymbol{\gamma}+t\boldsymbol{\bar{\gamma}}$ with respect to its first marginal $\beta_1=(1-t)\gamma_1+t\bar{\gamma}_1$, we get that 
\begin{align} \label{F-L convex 1}
\int_{X_1}C\left(x_1,\beta_{x_1}\right)d\beta_1&\leq (1-t)\int_{X_1}C(x_1, \gamma_{x_1})d\gamma_1+t\int_{X_1}C(x_1,\bar{\gamma}_{x_1}
)d\bar{\gamma}_1<\infty.
\end{align}		
This means that $M$ is a convex set.

		For every $\boldsymbol{\gamma}\in M$, applying Lemma \ref{L-inequality of F} we obtain that \begin{align*}
			&\cE_C(\boldsymbol{\gamma} \vert \mu_1,\mu_2)\\
			&=\sup_{(\varphi_1,\varphi_2)\in B}\left\{\sum_{i=1}^2\int_{X_i}F_i^\circ(\varphi_i)d\mu_i+\int_{X_1}\left(C(x_1,\gamma_{x_1})- \varphi_1(x_1)\right)d\gamma_1-\int_{X_2}\varphi_2(x_2)d\gamma_2\right\}\\
			&=\sup_{(\varphi_1,\varphi_2)\in B}\left\{\sum_{i=1}^2\int_{X_i}F_i^\circ(\varphi_i)d\mu_i+\int_{X_1}\left(C(x_1,\gamma_{x_1})-\varphi_1(x_1)-\gamma_{x_1}(\varphi_2)\right)d\gamma_1\right\}.
		\end{align*}
		We now define the function $L$ on $M\times B$ by $$L(\boldsymbol{\gamma},\varphi):=\sum_{i=1}^2\int_{X_i}F_i^\circ(\varphi_i)d\mu_i+\int_{X_1}\left(C(x_1,\gamma_{x_1})-\varphi_1(x_1)-\gamma_{x_1}(\varphi_2)\right)d\gamma_1, $$
		for every $\boldsymbol{\gamma}\in M, \varphi=(\varphi_1,\varphi_2)\in B$.
		This yields that $$\cE(\mu_1,\mu_2)=\inf_{\boldsymbol{\gamma}\in \cM(X_1\times X_2)}\cE_C(\boldsymbol{\gamma} \vert \mu_1,\mu_2)=\inf_{\boldsymbol{\gamma}\in M}\sup_{\varphi\in B}L(\gamma,\varphi).$$
		On the other hand, for every $\varphi=(\varphi_1,\varphi_2)\in B$ it is not difficult to check that $$\inf_{\boldsymbol{\gamma}\in M}\int_{X_1}\left(C(x_1,\gamma_{x_1})-\varphi_1(x_1)-\gamma_{x_1}(\varphi_2)\right)d\gamma_1(x_1)=\left\lbrace\begin{array}{ll} 0 &  \text{ if } \varphi\in \Lambda, \\ -\infty  
			& \text{ otherwise. }\end{array}\right.$$
		Therefore, we obtain that $$\sup_{\varphi\in B}\inf_{\boldsymbol{\gamma}\in M}L(\boldsymbol{\gamma},\varphi)=\sup_{(\varphi_1,\varphi_2)\in \Lambda}\sum_{i=1}^2\int_{X_i}F_i^\circ(\varphi_i)d\mu_i.$$
		Hence, we need to prove that $$\inf_{\boldsymbol{\gamma}\in M}\sup_{\varphi\in B}L(\boldsymbol{\gamma},\varphi)=\sup_{\varphi\in B}\inf_{\boldsymbol{\gamma}\in M}L(\boldsymbol{\gamma},\varphi).$$

		As $\varphi_i$ is continuous for $i=1,2$, and $C$ is lower semi-continuous one has $L(\cdot,\varphi)$ is lower semi-continuous on $M$. Since $F_i^\circ$ is concave, we get that $L$ is concave in $B$. Moreover, for any $\boldsymbol{\gamma},\boldsymbol{\bar{\gamma}}\in M$ and $t\in [0,1]$, putting $\boldsymbol{\beta}:=(1-t)\boldsymbol{\gamma}+t\boldsymbol{\bar{\gamma}}$, one has \begin{align*}
        &\int_{X_1}\beta_{x_1}(\varphi_2)d\beta_1(x_1)\notag \\
        = &(1-t)\int_{X_1}\dfrac{d\gamma_1}{d((1-t)\gamma_1+t\bar{\gamma}_1)}\gamma_{x_1}(\varphi_2)d((1-t)\gamma_1+t\bar{\gamma}_1)(x_1)\notag\\
        &+t\int_{X_1}\dfrac{d\bar{\gamma}_1}{d((1-t)\gamma_1+t\bar{\gamma}_1)}\bar{\gamma}_{x_1}(\varphi_2)d((1-t)\gamma_1+t\bar{\gamma}_1)(x_1)\notag\\
        =& (1-t)\int_{X_1}\gamma_{x_1}(\varphi_2)d\gamma_1(x_1)+t\int_{X_1}\bar{\gamma}_{x_1}(\varphi_2)d\bar{\gamma}_1(x_1).   	
		\end{align*}
Combining with \eqref{F-L convex 1} we obtain that the function $L$ is convex on $M$.

		Next, from the coercive condition $(F_1)'_\infty+(F_2)'_\infty+\inf C>0$, we can find constant functions $\overline{\varphi}_i\in (-(F_i)'_\infty,+\infty)$ for $i=1,2$ such that $\inf C-\overline{\varphi}_1-\overline{\varphi}_2>0$. Then let $\overline{\varphi}=(\overline{\varphi}_1,\overline{\varphi}_2)$, for every $\boldsymbol{\gamma}\in M$, one has $$L(\boldsymbol{\gamma},\overline{\varphi})=\sum_{i=1}^2F_i^\circ(\overline{\varphi}_i) \vert \mu_i \vert +\int_{X_1}\left(C(x_1,\gamma_{x_1})-\inf C\right)d\gamma_1+(\inf C-\overline{\varphi}_1-\overline{\varphi}_2)\gamma_1(X_1).$$
		This implies that for large enough $K>0$ we get that $D:=\{\boldsymbol{\gamma}\in M \vert L(\boldsymbol{\gamma},\overline{\varphi})\leq K\}$ is bounded. As $X_i$ is compact one has $D$ is also equally tight. Hence, using Prokhorov's Theorem we obtain that $D$ is relatively compact under the weak topology. Observe that as $L(\cdot,\overline{\varphi})$ is lower semi-continuous one has $D$ is closed. Therefore, by \cite[Theorem 2.4]{LMS18} we get the result. 
	   \end{proof}
	
When $X_1,X_2$ may not be compact, to obtain our duality formula as in Theorem \ref{dual theorem} we need new ideas which are different from \cite{LMS18}. Our proof of Theorem 1 relies on the fact that the functional $\ET$, defined in \eqref{D-functional ET}, is convex and positively homogenous and lower semi-continuous, and is thus the support function of a convex set. This fact is established in Lemma \ref{L-lower semicontinuous of ET} Lemma \ref{L-Lamda_{ET}^<=Lamda_R}. 
	\begin{lemma}\label{L-easy part of the dual formula} Let $X_1,X_2$ be Polish metric spaces and assume that $(F_1)'_\infty+(F_2)'_\infty+\inf C>0$. Then for every $\mu_i\in \cM(X_i), i=1,2$ we have 
		$$\cE_C(\mu_1,\mu_2)\geq\sup_{(\varphi_1,\varphi_2)\in \Lambda}\sum_{i=1}^2\int_{X_i}F_i^\circ(\varphi_i)d\mu_i.$$
In particular, $\cE(\cdot \vert \mu_1,\mu_2)$ is bounded from below. 		
	\end{lemma}
	\begin{proof}
	For every $(\varphi_1,\varphi_2)\in \Lambda$ and $\boldsymbol{\gamma}\in\cM(X_1\times X_2)$, applying Lemma \ref{L-inequality of F} we get that 
		\begin{eqnarray*}
			\cE_C(\boldsymbol{\gamma} \vert \mu_1,\mu_2)&=&\cF_1(\gamma_1 \vert \mu_1)+\cF_2(\gamma_2 \vert \mu_2)+\int_{X_1}C(x_1,\gamma_{x_1})d\gamma_1(x_1)\\
			&\geq & \cF_1(\gamma_1 \vert \mu_1)+\cF_2(\gamma_2 \vert \mu_2)+\int_{X_1}(\varphi_1(x_1)+\gamma_{x_1}(\varphi_2))d\gamma_1(x_1)\\
			&=& \cF_1(\gamma_1 \vert \mu_1)+\int_{X_1}\varphi_1d\gamma_1+\cF_2(\gamma_2 \vert \mu_2)+\int_{X_1}\int_{X_2}\varphi_2(x_2)d\gamma_{x_1}(x_2)d\gamma_1(x_1)\\
			&=& \cF_1(\gamma_1 \vert \mu_1)+\int_{X_1}\varphi_1d\gamma_1+\cF_2(\gamma_2 \vert \mu_2)+\int_{X_2}\varphi_2d\gamma_2\\
			&\geq & \sum_{i=1}^2\int_{X_i}F_i^\circ(\varphi_i)d\mu_i.
		\end{eqnarray*}

Next, since the condition $(F_1)'_\infty+(F_2)'_\infty+\inf C>0$, by the same way in the proof of Lemma \ref{L-the set lamda is nonempty} we can find constant functions $\overline{\varphi}_1$ and $\overline{\varphi}_2$ such that $(\overline{\varphi}_1,\overline{\varphi}_2)\in\Lambda$. Therefore,
$$\cE_C(\boldsymbol{\gamma} \vert \mu_1,\mu_2)\geq \sum_{i=1}^2F_i^\circ(\overline{\varphi}_i) \vert \mu_i \vert > -\infty.$$
	   \end{proof}
	
	Recall that given a metric space $X$ we denote by $(C_b(X))^*$ the dual space of the normed space $(C_b(X),\|\cdot\|_\infty)$. 	We recall the functional $\ET:(C_b(X_1))^*\times (C_b(X_2))^*\to [-\infty,+\infty]$ we defined in \eqref{D-functional ET}.
\begin{align*}
\ET(T_1,T_2):=\left\lbrace\begin{array}{ll}
		\cE_C(\mu_1,\mu_2) & \text{ if } (T_1,T_2)=(T_{\mu_1}, T_{\mu_2}) 
		,\\
		+\infty &\text{ otherwise.}
	\end{array}\right.
\end{align*}

	Note that if $(F_1)'_\infty+(F_2)'_\infty+\inf C>0$ then by Lemma \ref{L-easy part of the dual formula} we always have $\ET(\mu_1,\mu_2)\neq -\infty$ for every $(\mu_1,\mu_2)\in \cM(X_1)\times\cM(X_2)$. 
	We define 
	\begin{align*}
	\Lambda_{ET} :=\bigg\{&(\varphi_1,\varphi_2)\in C_b(X_1)\times C_b(X_2):\sum_{i=1}^2\int_{X_i}\varphi_id\mu_i\leq \ET(\mu_1,\mu_2),\\
	& \mbox{ for every } (\mu_1,\mu_2)\in \cM(X_1)\times\cM(X_2)\bigg\}, 
	\end{align*}
	$$\Lambda_{ET}^<:=\{(\varphi_1,\varphi_2)\in \Lambda_{ET} \vert \sup_{x_i\in X_i}\varphi_i(x_i)<F_i(0),i=1,2\}.$$
	\begin{lemma}\label{L-Lamda_ET<F(0)}
		Let $X_1,X_2$ be Polish metric spaces. If $\Lambda_{ET}^<$ is a nonempty set then for every $\mu_i\in \cM(X_i),i=1,2$ one has $$\sup_{(\varphi_1,\varphi_2)\in \Lambda_{ET}}\sum_{i=1}^2\int_{X_i}\varphi_id\mu_i=\sup_{(\varphi_1,\varphi_2)\in \Lambda_{ET}^<}\sum_{i=1}^2\int_{X_i}\varphi_id\mu_i$$
	\end{lemma}
	\begin{proof}
		It is clear that we only need to show that $$\sup_{(\varphi_1,\varphi_2)\in \Lambda_{ET}}\sum_{i=1}^2\int_{X_i}\varphi_id\mu_i\leq\sup_{(\varphi_1,\varphi_2)\in \Lambda_{ET}^<}\sum_{i=1}^2\int_{X_i}\varphi_id\mu_i.$$ For every $\varepsilon>0$ there exists $(\phi_1^\varepsilon,\phi_2^\varepsilon)\in \Lambda_{ET}$ such that $$\sum_{i=1}^2\int_{X_i}\phi_i^\varepsilon d\mu_i\geq \sup_{(\varphi_1,\varphi_2)\in \Lambda_{ET}}\sum_{i=1}^2\int_{X_i}\varphi_id\mu_i-\varepsilon/2.$$
If $ \vert \mu_1 \vert = \vert \mu_2 \vert =0$, we are done. Otherwise, for each $i\in\{1,2\}$, set $\overline{\phi}_i^\varepsilon:=\phi^\varepsilon_i-\varepsilon/(2( \vert \mu_1 \vert + \vert \mu_2 \vert ))$.

		Moreover, denote by $\eta$ the null measure on $X_1\times X_2$. As $(\phi_1^\varepsilon,\phi_2^\varepsilon)\in \Lambda_{ET}$, for every $(\mu_1,\mu_2)\in \cM(X_1)\times\cM(X_2)$ one has
		 $$\sum_{i=1}^2\int_{X_i}\phi_i^\varepsilon d\mu_i\leq \ET(\mu_1,\mu_2)\leq \cE_C(\eta \vert \mu_1,\mu_2)=F_1(0)\mu_1(X_1)+F_2(0)\mu_2(X_2).$$
		For any $x_1\in X_1$ set $\mu_1:=\delta_{x_1}$ and $\mu_2$ is the null measure on $X_2$ we get that $\phi_1^\varepsilon(x_1)\leq F_1(0)$. Similarly, we also have $\phi_2^\varepsilon(x_2)\leq F_2(0)$ for every $x_2\in X_2$. Therefore, $$\sup_{x_i\in X_i}\overline{\phi}_i^\varepsilon (x_i) =\sup_{x_i\in X_i}\phi_i^\varepsilon (x_i)-\frac{\varepsilon}{2( \vert \mu_1 \vert + \vert \mu_2 \vert )}<F_i(0),\,i=1,2.$$ Thus,$(\overline{\phi}_1^\varepsilon,\overline{\phi}_1^\varepsilon)\in \Lambda_{ET}^<$. Hence, we obtain that \begin{align*}
			\sup_{(\varphi_1,\varphi_2)\in \Lambda_{ET}^<}\sum_{i=1}^2\int_{X_i}\varphi_id\mu_i &\geq \sum_{i=1}^2\int_{X_i}\overline{\phi}_{i}^\varepsilon d\mu_i\\
			&=\sum_{i=1}^2\int_{X_i}\phi_{i}^\varepsilon d\mu_i-\varepsilon/2\\
			&\geq \sup_{(\varphi_1,\varphi_2)\in \Lambda_{ET}}\sum_{i=1}^2\int_{X_i}\varphi_id\mu_i-\varepsilon.
		\end{align*}
		So that the proof is complete.
	   \end{proof}
	
	\begin{lemma}\label{L-lower semicontinuous of ET}
		Let $X_1,X_2$ be Polish spaces and $\mu_i\in \cM(X_i),i=1,2$. Assume that $F_i$ is superlinear for $i=1,2$. For each $i\in \{1,2\}$, let $(\mu_{i}^n)_n\subset\cM(X_i)$  such that $\mu_{i}^n$ converges to $\mu_i$ in the weak topology then $$\liminf_{n\rightarrow \infty}\cE_C(\mu_{1}^n,\mu_{2}^n)\geq \cE_C(\mu_1,\mu_2).$$ 
	\end{lemma}
	\begin{proof} If $\liminf_{n\rightarrow \infty}\cE_C(\mu_{1}^n,\mu_{2}^n)=+\infty$, we are done. Otherwise, we can assume that $\cE_C(\mu_{1}^n,\mu_{2}^n)<M<\infty$ for every $n\in\mathbb{N}$. For each $n\in \mathbb{N}$, using Theorem \ref{T-existence of minimizers} let $\boldsymbol{\gamma^n}\in \cM(X_1\times X_2)$ such that $\cE_C(\mu_{1}^n,\mu_{2}^n)=\cE(\boldsymbol{\gamma^n} \vert \mu_{1}^n,\mu_{2}^n)$. As $\mu_{i}^n$ converges to $\mu_i$ one has $(\mu_{i}^n)_n$ is bounded and equally tight for $i=1,2$. Moreover, observe that for $i\in\{1,2\}$ we have $\cF(\gamma_{i}^n \vert \mu_{i}^n)\leq \cE_C(\mu_{1}^n,\mu_{2}^n)<M$ for every $n\in \mathbb{N}$. Hence, applying \cite[Proposition 2.10]{LMS18} we get that $(\gamma_i^n)_n$ is equally tight and bounded for $i=1,2$. By \cite[Lemma 5.2.2]{AGS} one gets that $(\boldsymbol{\gamma^n})_n$ is also equally tight and bounded. Therefore, by Prokhorov's Theorem, passing to a subsequence we can assume that $\boldsymbol{\gamma^n}\rightarrow\boldsymbol{\gamma}$ as $n\rightarrow\infty$ in the weak topology for some $\boldsymbol{\gamma}\in \cM(X_1\times X_2)$. From Lemma \ref{L-inequality of F} we get that the function $\cF$ is lower semi-continuous. This implies that $$\liminf_{n\rightarrow\infty}\sum_{i=1}^2\cF_i(\gamma^n_i \vert \mu_{i}^n)\geq \sum_{i=1}^2\cF_i(\gamma_i \vert \mu_i).$$
		Next, applying Lemma \ref{L-lower semicontinuous of integration C} we also obtain that $$\liminf_{n\rightarrow\infty}\int_{X_1}C(x_1,\gamma_{x_1}^n)d\gamma_{1}^n(x_1)\geq \int_{X_1}C(x_1,\gamma_{x_1})d\gamma_1(x_1).$$
		Therefore, we get the result.
	   \end{proof}
In our previous manuscript, we added a technical condition called (BM), playing a crucial role in our work. It turns out that property (BM) is always true, as shown in the following lemma. The statement of this lemma and its proof are suggested by one of the referees. We thank her/him for the suggestion, which has helped us completely remove this technical condition and has significantly improved our results.	
\begin{lemma}\label{L-BM}
Let $X$ be a Polish metric space and $F:[0;+\infty)\to [0;+\infty)$ be a convex lower semi-continuous function. Let $R:[0,\infty)\to [0,\infty]$ be the reverse density function of $F$, i.e. $R(r)=rF(1/r)$ if $r>0$ and $R(0)=F_\infty'$. Then for every $\psi\in C_b(X)$ satisfying that $\psi(x)\in (-\mathrm{aff } F_\infty, F(0))$ for all $x\in X$ with 
\begin{align*}
\mathrm{aff } F_\infty=\begin{cases}
+\infty &\mbox{ if } F'_\infty=+\infty,\\
\lim_{u\to \infty}F_\infty' u-F(u) &\mbox{ otherwise},
\end{cases}
\end{align*} there exists a Borel bounded function $s:X\rightarrow (0,\infty)$ such that \begin{align*}
			R(s(x))+R^*(\psi(x))=s(x)\psi(x),\text{ for every }x\in X. 	
			\end{align*}
			In particular, if $F$ is superlinear, i.e. $F'_\infty=+\infty$ then for every $\psi\in C_b(X)$ such that $\sup_{x\in X}\psi(x)<F(0)$, there exists a Borel bounded function $s:X\rightarrow (0,\infty)$ such that \begin{align*}
			R(s(x))+R^*(\psi(x))=s(x)\psi(x),\text{ for every }x\in X. 	
			\end{align*}
\end{lemma}
\begin{proof}
First, we extend the function $R$ by $\widetilde{R}:\mathbb{R} \to (-\infty;+\infty]$ as 
\begin{align*}
\widetilde{R}(r)=\begin{cases}
R(r) & \text{ if }r\geq 0,\\
+\infty & \text{ if }r < 0.
\end{cases}
\end{align*}
Then $R^*$ is the conjugate of $\widetilde{R}$. Observe that $\widetilde{R}$ is  convex and lower semi-continuous, thus applying \cite[Proposition 3.1, page 14 and Proposition 4.1, page 18]{ET} one gets that $(R^*)^*=\widetilde{R}$. Hence, by \cite[(2.17)]{LMS18}, for every $t\in D(R^*)$, if $s\in \partial R^*(t)$ then we have $s\in D((R^*)^*)$ and \begin{align*}
R^*(t)+R(s)=R^*(t)+\widetilde{R}(s)=R^*(t)+(R^*)^*(s)=st,
\end{align*}
where $\partial R^*$ is the subdifferential of $R^*$ at $t$.

Recall that $\mathring{D}(R^*)$ is the interior of the domain of $R^*$. As $R^*$ is convex, for any $t_0\in \mathring{D}(R^*)$ we get that the left derivative of $R^*$ at $t_0$, $D_{-}R^*(t_0)=\lim_{t\to t^{-}_0}\dfrac{R^*(t)-R^*(t_0)}{t-t_0}$ exists and $D_{-}R^*(t_0)\in \partial R^*(t_0)$. Furthermore, as $R^*$ is continuous on $\mathring{D}(R^*)$, we get that the function $t\mapsto D_{-}R^{*}(t)$ is measurable. So if $\psi$ is some bounded continuous function taking values in $\mathring{D}(R^*)$, the equality 
\begin{align*}
			R(s(x))+R^*(\psi(x))=s(x)\psi(x),\text{ for every }x\in X,	
			\end{align*}
holds with $s(x)=D_{-}R^*(\psi(x))$, which is measurable, as a composition of measurable maps.

Now, let $\psi\in C_b(X)$ with $\psi(x)\in (-\mathrm{aff } F_\infty, F(0))$ for all $x\in X$, then $\psi(x)\in \mathring{D}(R^*)$ for all $x\in X$, since $\mathring{D}(R^*)=(-\infty, F(0))$ (see \cite[the first line in page 992]{LMS18}). Since $R^*$ is strictly increasing on $(-\mathrm{aff } F_\infty, F(0))$ (see \cite[(2.33)]{LMS18}), one gets $s(x)=D_{-}R^*(\psi(x))>0$ for all $x\in X$.

Finally, since the map $t\mapsto D_{-}R^{*}(t)$ is increasing and $\psi$ is upper bounded, the function $s$ is also upper bounded, which completes the proof.
\end{proof}	
For the convenience, we will write $\ET(\mu_1,\mu_2)$ for $\ET(T_{\mu_1}, T_{\mu_2})$	for every $(\mu_1,\mu_2)\in \cM(X_1)\times \cM(X_2)$.
	\begin{lemma}
	\label{L-Lamda_{ET}^<=Lamda_R} Let $X_1,X_2$ be Polish metric spaces. Suppose that $F_1,F_2$ are superlinear. Then
\begin{enumerate}
\item the functional $\ET:(C_b(X_1))^*\times(C_b(X_2))^*\to (-\infty,+\infty]$ is convex and positively one homogeneous, i.e. $\ET(\lambda T_1,\lambda T_2)=\lambda \ET(T_1,T_2)$ for every $\lambda\geq 0, T_1\in (C_b(X_1))^*,T_2\in (C_b(X_2))^*$;
\item $\Lambda_R\subset\Lambda_{ET}^<$;
\item $\Lambda_{ET}^<=\Lambda_R$.
\end{enumerate}	
	\end{lemma}
	\begin{proof}
1) By the construction of $\ET$, it is clear that $\ET(0,0)=0$ and $\ET(\lambda T_1,\lambda T_2)=\lambda \ET(T_1,T_2)$ for every $\lambda\geq 0$, $(T_1,T_2)\notin \cM(X_1)\times \cM(X_2)$ (here we use the convention that $0\cdot (+\infty)=0)$. Therefore, to check
$\ET$ is positively one homogeneous we only need to check $\ET(\lambda T_{\mu_1},\lambda T_{\mu_2})=\lambda \ET(T_{\mu_1},T_{\mu_2})$ for every $\lambda>0$, $(\mu_1,\mu_2)\in \cM(X_1)\times \cM(X_2).$
Given $\boldsymbol{\gamma}\in \cM(X_1\times X_2)$ and $\lambda>0$ then its disintegration $(\gamma_{x_1})_{x_1\in X_1}$ with respect to the first marginal $\gamma_1$ is also the disintegration of $\lambda \boldsymbol{\gamma}$ with respect to its first marginal $\lambda \gamma_1$.  From Lemma \ref{L-inequality of F}, one has that $\cF_i(\lambda \gamma_i\vert\lambda \mu_i)=\lambda \cF_i(\gamma_i\vert\mu_i)$ for $i=1,2$. Hence
\begin{align*}
\ET(\lambda T_{\mu_1},\lambda T_{\mu_2})&=\ET(T_{\lambda \mu_1}, T_{\lambda\mu_2})=\cE_C(\lambda\mu_1,\lambda\mu_2)\\
&=\inf\{\cE_C(\boldsymbol{\gamma}\vert\lambda\mu_1,\lambda\mu_2): \boldsymbol{\gamma}\in \cM(X_1\times X_2)\}\\
&=\inf\{\cE_C(\lambda\boldsymbol{\gamma}\vert\lambda\mu_1,\lambda\mu_2): \boldsymbol{\gamma}\in \cM(X_1\times X_2)\}\\
&=\inf\big\{\sum_{i=1}^2\cF_i(\lambda\gamma_i\vert\lambda\mu_i)+\lambda\int_{X_1}C(x_1,\gamma_{x_1})d\gamma_1(x_1): \boldsymbol{\gamma}\in \cM(X_1\times X_2)\}\\
&=\lambda\inf\big\{\sum_{i=1}^2\cF_i(\gamma_i\vert\mu_i)+\int_{X_1}C(x_1,\gamma_{x_1})d\gamma_1(x_1): \boldsymbol{\gamma}\in \cM(X_1\times X_2)\big\}\\
&=\lambda \cE_C(\mu_1,\mu_2)=\lambda \ET(T_{\mu_1},T_{\mu_2}).
\end{align*}
 Since the homogeneity property of $\ET$, to show that $\ET$ is convex, we only need to check that $$\ET(\mu_1,\mu_2)+\ET(\nu_1,\nu_2)\geq \ET(\mu_1+\nu_1,\mu_2+\nu_2) \mbox{ for every } \mu_i,\nu_i\in \cM(X_i),i=1,2.$$ We will consider $(\mu_1,\mu_2),(\nu_1,\nu_2)\in \cM(X_1)\times \cM(X_2)$ such that $ \cE_C(\mu_1,\mu_2)<\infty$ and $\cE_C(\nu_1,\nu_2)<\infty $ (the other cases are trivial). From Theorem \ref{T-existence of minimizers}, let $\boldsymbol{\gamma},\overline{\boldsymbol{\gamma}}\in \cM(X_1\times X_2)$ such that $\ET(\mu_1,\mu_2)=\cE_C(\boldsymbol{\gamma} \vert \mu_1,\mu_2)$ and $\ET(\nu_1,\nu_2)=\cE_C(\overline{\boldsymbol{\gamma}} \vert \nu_1,\nu_2)$. 
				
		As $\bigg((d\gamma_1/d(\gamma_1+\overline{\gamma}_1))\gamma_{x_1}+(d\overline{\gamma}_1/d(\gamma_1+\overline{\gamma}_1))\overline{\gamma}_{x_1}\bigg)_{x_1\in X_1}$ is the disintegration of $\boldsymbol{\gamma}+\overline{\boldsymbol{\gamma}}$ with respect to $\gamma_1+\overline{\gamma}_1$ and $C(x_1,\cdot)$ is convex on $\cP(X_2)$ for every $x_1\in X_1$, we obtain that 
		$$\int_{X_1}C(x_1,\gamma_{x_1})d\gamma_1+\int_{X_1}C(x_1,\overline{\gamma}_{x_1})d\overline{\gamma}_1\geq \int_{X_1}C(x_1, (\gamma+\overline{\gamma})_{x_1})d(\gamma_1+\overline{\gamma}_1).$$ 		
		This implies that 
		\begin{align*}
			\ET(\mu_1,\mu_2)+\ET(\nu_1,\nu_2)&=\sum_{i=1}^2\left(\cF_i(\gamma_i \vert \mu_i)+\cF_i(\overline{\gamma}_i \vert \nu_i)\right)\\
			&+\int_{X_1}C(x_1,\gamma_{x_1})d\gamma_1+\int_{X_1}C(x_1,\overline{\gamma}_{x_1})d\overline{\gamma}_1\\
			&\geq \sum_{i=1}^2\cF_i(\gamma_i+\overline{\gamma}_i \vert 
			\mu_i+\nu_i)+\int_{X_1}C(x_1, (\gamma+\overline{\gamma})_{x_1})d(\gamma_1+\overline{\gamma}_1)\\
			&\geq \ET(\mu_1+\nu_1,\mu_2+\nu_2).
		\end{align*} 		
		Therefore, $\ET$ is convex.

2) Let any $(\varphi_1,\varphi_2)\in \Lambda_R$. Now we will prove that $\sum_{i=1}^2\int_{X_i}\varphi_id\mu_i\leq\ET(\mu_1,\mu_2),$ for every $(\mu_1,\mu_2)\in \cM(X_1)\times\cM(X_2)$.
If $(\mu_1,\mu_2)\in \cM(X_1)\times\cM(X_2)$ such that $\ET(\mu_1,\mu_2)=+\infty $ then it is clear. So we only consider the case $(\mu_1,\mu_2)\in \cM(X_1)\times\cM(X_2)$ such that $ \cE_C(\mu_1,\mu_2)<\infty $. From Theorem \ref{T-existence of minimizers}, let $\boldsymbol{\gamma}\in \cM(X_1\times X_2)$ such that $\ET(\mu_1,\mu_2)=\cE_C(\boldsymbol{\gamma} \vert \mu_1,\mu_2)$. Then we have that 
		\begin{align*}
			\cE_C(\boldsymbol{\gamma} \vert \mu_1,\mu_2)&=\sum_{i=1}^2\cF_i(\gamma_i \vert \mu_i)+\int_{X_1}C(x_1,\gamma_{x_1})d\gamma_1(x_1)\\
			&\geq \sum_{i=1}^2\cF_i(\gamma_i \vert \mu_i)+\int_{X_1}\left( R_1^*(\varphi_1(x_1))+\gamma_{x_1}(R_2^*(\varphi_2))\right)d\gamma_1(x_1)	\\
			&=\sum_{i=1}^2\cF_i(\gamma_i \vert \mu_i)+\int_{X_1}R_1^*(\varphi_1(x_1))d\gamma_1(x_1)+\int_{X_1}\int_{X_2}R_2^*(\varphi_2(x_2))d\gamma_{x_1}(x_2)d\gamma_1(x_1)\\
			&=\sum_{i=1}^2\cF_i(\gamma_i \vert \mu_i)+\int_{X_1}R_1^*(\varphi_1(x_1))d\gamma_1(x_1)+\int_{X_2}R_2^*(\varphi_2(x_2))d\gamma_2(x_2).
		\end{align*}
		 Applying Lemma \ref{L-inequality of F} we get that 
		$$\int_{X_i}\varphi_id\mu_i\leq \cF_i(\gamma_i \vert \mu_i)+\int_{X_i}R^*_i(\varphi_i)d\gamma_i.$$		
		Therefore, $$\sum_{i=1}^2\int_{X_i}\varphi_id\mu_i\leq \cE_C(\boldsymbol{\gamma} \vert \mu_1,\mu_2)=\ET(\mu_1,\mu_2).$$
		This implies that $(\varphi_1,\varphi_2)\in \Lambda^<_{ET}$. Hence, $\Lambda_R\subset \Lambda^<_{ET}$. This also shows that $\Lambda^<_{ET}$ is nonempty. 
		
3)	We now check that $\Lambda_{ET}^<=\Lambda_R$. 		
		
		Let any $(\varphi_1,\varphi_2)\in \Lambda_{ET}^<$. For every $\overline{x}_1\in X_1, p\in \cP(X_2), r>0$ we define $\mu_1:=\delta_{\overline{x}_1}$ and $\boldsymbol{\gamma}:=r\delta_{\overline{x}_1}\otimes p$ then for every $\mu_2\in \cM(X_2)$ one has \begin{align*}
			\varphi_1(\overline{x}_1)+\int_{X_2}\varphi_2(x_2)d\mu_2(x_2)&\leq \ET(\mu_1,\mu_2)\\
			&\leq \cE_C(\boldsymbol{\gamma} \vert \mu_1,\mu_2)\\
			&=F_1(r)+\cF(\gamma_2 \vert \mu_2)+rC(\overline{x}_1,p).
		\end{align*}
		This yields, \begin{align*}
			\frac{1}{r}\left(\varphi_1(\overline{x}_1)-F_1(r)\right)\leq C(\overline{x}_1,p)+\frac{1}{r}\left(\cF(\gamma_2 \vert \mu_2)-\int_{X_2}\varphi_2d\mu_2\right),\, \forall\mu_2\in \cM(X_2).
		\end{align*}
		Applying Lemma \ref{L-BM}, there exists a Borel bounded function $s:X_2\rightarrow (0,\infty)$ such that \begin{align}\label{F-BM1}
R(s(x))+R^*(\varphi_2(x))=s(x)\varphi_2(x),\text{ for every }x\in X.
\end{align}
		 Next, put $\mu_2:=s\gamma_2$. As $s$ is Borel bounded function one has $\mu_2\in \cM(X_2)$. We will check that $\gamma_2$ is absolutely continuous w.r.t $\mu_2$. For every Borel subset $A$ of $X_2$ such that $\mu_2(A)=0$ one has $\int_As(x)d\gamma_2=0$. Notice that $s(x)>0$ for every $x\in A$, hence $\gamma_2(A)=0$. So $\gamma_2$ is absolutely continuous w.r.t $\mu_2$.\\
		As $\varphi_2$ is bounded and $\sup_{x_2\in X_2}\varphi_2(x_2)<F_2(0)$, applying \eqref{F-interior of R*} we get that $R^*_2(\varphi_2)$ is bounded by $ R_2^*(\inf\varphi_2), R_2^*(\sup\varphi_2)\in  \mathbb{R}  $. Thus, from (\ref{F-BM1}) we obtain that $R_2(s)$ is also bounded. Hence, by \eqref{F-relation between entropy and its reverse} one has $$\cF_2(\gamma_2 \vert \mu_2)=\cR(\mu_2 \vert \gamma_2)=\int_{X_2}R_2(s(x_2))d\gamma_2(x_2)<\infty.$$		
		Therefore, applying Lemma \ref{L-equality of F for (BM)} we obtain that $$\cF_2(\gamma_2 \vert \mu_2)-\int_{X_2}\varphi_2d\mu_2=-\int_{X_2}R_2^*(\varphi_2)d\gamma_2.$$
		Hence, for every $\overline{x}_1\in X_1,p\in \cP(X_2)$ and $r>0$ we get that \begin{align*}
			\frac{1}{r}\left(\varphi_1(\overline{x}_1)-F_1(r)\right)\leq C(\overline{x}_1,p)-\frac{1}{r}\int_{X_2}R_2^*(\varphi_2)d\gamma_2.
		\end{align*}
		Furthermore, observe that $\gamma_2=rp$ we obtain $$R^*_1(\varphi_1(x_1))=\sup_{r>0}\left(\varphi_1(x_1)-F_1(r)\right)/r\leq C(x_1,p)-p\left(R_2^*(\varphi_2)\right),\,\forall x_1\in X_1,p\in \cP(X_2).$$
		This implies that $\Lambda_{ET}^<\subset \Lambda_R$ and thus we get the result.
	   \end{proof}
	
	\begin{lemma}\label{L-dual theorem}
		Let $X_1,X_2$ be Polish metric spaces. Assume that $F_i$ is superlinear for $i=1,2$. Then for every $\mu_i\in \cM(X_i),i=1,2$ we have that
		\begin{align*}
		\cE_C(\mu_1,\mu_2)&=\sup_{(\varphi_1,\varphi_2)\in\Lambda_{ET}}\sum_{i=1}^2\int_{X_i}\varphi_id\mu_i =\sup_{(\varphi_1,\varphi_2)\in\Lambda_{ET}^<}\sum_{i=1}^2\int_{X_i}\varphi_id\mu_i\\
		&=	 \sup_{(\varphi_1,\varphi_2)\in \Lambda_R}\sum_{i=1}^2\int_{X_i}\varphi_id\mu_i.
		\end{align*}
	\end{lemma}
	\begin{proof}
		Since the one homogeneity of $\ET$ (see Lemma \ref{L-Lamda_{ET}^<=Lamda_R}), it is not difficult to check that $$\ET^*(\varphi_1,\varphi_2)=\left\lbrace\begin{array}{ll}
			0 &\text{ if } (\varphi_1,\varphi_2)\in \Lambda_{ET}^0,\\
			+\infty &\text{ otherwise,}
		\end{array}\right.$$
		where $\Lambda_{ET}^0:=\{(\varphi_1,\varphi_2)\in C_b(X_1)\times C_b(X_2) \vert (\varphi_1,\varphi_2)\in\Lambda_{ET}\}$. 

		Moreover, by Lemmas \ref{L-lower semicontinuous of ET} and \ref{L-Lamda_{ET}^<=Lamda_R} one has $\ET$ is convex and lower semi-continuous under the weak topopology. Hence, by \cite[Proposition 3.1, page 14 and Proposition 4.1, page 18]{ET} we get that $\left(\ET^*\right)^*=\ET$. Therefore,
		\begin{align*}\sup_{(\varphi_1,\varphi_2)\in\Lambda_{ET}}\sum_{i=1}^2\int_{X_i}\varphi_id\mu_i & \leq \ET(\mu_1,\mu_2)\\
			&= (\ET^*)^*(\mu_1,\mu_2)\\
			&=\sup_{(\varphi_1,\varphi_2)\in C_b(X_1)\times C_b(X_2)}\left\{\sum_{i=1}^2\int_{X_i}\varphi_id\mu_i-\ET^*(\varphi_1,\varphi_2)\right\}\\
			&= \sup_{(\varphi_1,\varphi_2)\in\Lambda_{ET}^0}\sum_{i=1}^2\int_{X_i}\varphi_id\mu_i\\
			&\leq \sup_{(\varphi_1,\varphi_2)\in\Lambda_{ET}}\sum_{i=1}^2\int_{X_i}\varphi_id\mu_i.
		\end{align*}		
		This implies that $\ET(\mu_1,\mu_2)=\sup_{(\varphi_1,\varphi_2)\in\Lambda_{ET}}\sum_{i=1}^2\int_{X_i}\varphi_id\mu_i.$ Thus, using Lemma \ref{L-Lamda_ET<F(0)} and Lemma \ref{L-Lamda_{ET}^<=Lamda_R} we obtain that \begin{align*}
			\cE_C(\mu_1,\mu_2)&=\ET(\mu_1,\mu_2) =\sup_{(\varphi_1,\varphi_2)\in\Lambda_{ET}}\sum_{i=1}^2\int_{X_i}\varphi_id\mu_i\\
			&=\sup_{(\varphi_1,\varphi_2)\in\Lambda_{ET}^<}\sum_{i=1}^2\int_{X_i}\varphi_id\mu_i\\
			&=\sup_{(\varphi_1,\varphi_2)\in\Lambda_R}\sum_{i=1}^2\int_{X_i}\varphi_id\mu_i.
		\end{align*}
	   \end{proof}
	\begin{lemma}\label{L-dual formula for Lamda set} Assume that $F_1, F_2$ are superlinear. Then
		\[\cE_C(\mu_1,\mu_2)= \sup_{(\varphi_1,\varphi_2)\in \Lambda_{R}}\sum_{i=1}^2\int_{X_i}\varphi_id\mu_i=\sup_{(\varphi_1,\varphi_2)\in \Lambda}\sum_{i=1}^2\int_{X_i}F_i^\circ(\varphi_i)d\mu_i. \]
	\end{lemma}
	\begin{proof}
		For $ (\varphi_1,\varphi_2)\in \Lambda_{R}$ and $i\in \{1,2\}$, we define $ \overline{\varphi}_i=R_i^*(\varphi_i).$ Because $ F_i $ is superlinear then $ \mathring{D}(F_i^\circ)=  \mathbb{R} . $ Because $ \varphi _i$ is bounded below by some number $ M_i<F_i(0) $ so $ \overline{\varphi}_i$ is bounded below by $ R_i^*(M)>-\infty. $ We have $ \overline{\varphi}_i $ is bounded above by $ R_i^*(\sup_{x_i\in X_i}\varphi_i(x_i)). $ To confirm $ (\overline{\varphi}_1,\overline{\varphi}_2)\in \Lambda $ we see that \begin{align*}
			\overline{\varphi}_1(x_1)+p(\overline{\varphi}_2)=R_i^*(\varphi_1(x_1))+p(R_i^*(\varphi_2))\le C(x_1,p),
		\end{align*}
		for every $x_1\in X_1$ and $p\in \cP(X_2)$. As $F_i^\circ(\overline{\varphi}_i)= F_i^\circ(R_i^*(\varphi_i))\ge\varphi_i $ (\cite[(2.31)]{LMS18}) one has $$\sum_{i=1}^2\int_{X_i}\varphi_id\mu_i\leq \sum_{i=1}^2\int_{X_i}F_i^\circ(\overline{\varphi}_i)d\mu_i.$$
		Thus, by Lemma \ref{L-dual theorem} and Lemma \ref{L-easy part of the dual formula} we get that 
		\[ \cE_C(\mu_1,\mu_2)=\sup_{(\phi_1,\phi_2)\in \Lambda_{R}}\sum_{i=1}^2\int_{X_i}\phi_id\mu_i\le\sup_{(\phi_1,\phi_2)\in \Lambda}\sum_{i=1}^2\int_{X_i}F_i^\circ(\phi_i)d\mu_i\le\cE_C(\mu_1,\mu_2) \] and the equalities happen.
	   \end{proof}
	\begin{lemma}\label{reduce to RC}
		We define $ R_C\varphi(x_1):=\inf_{p\in\cP(X_2)}\{C(x_1,p)+p(\varphi) \} $ for every $ x_1\in X_1 $ and $ \varphi\in C_b(X_2,\mathring{D}(F_2^\circ)). $
		Assume that $F_1$ and $F_2$ are superlinear. Then for every $\mu_i\in \cM(X_i),i=1,2$ one has 
\begin{align*}
\cE_C(\mu_1,\mu_2)&= \sup_{\varphi\in C_b(X_2)}\int_{X_1}F_1^\circ(R_C\varphi)d\mu_1+\int_{X_2}F_2^\circ(-\varphi)d\mu_2\\
&= \sup_{(\varphi_1,\varphi_2)\in \Lambda}\sum_{i=1}^2\int_{X_i}F_i^\circ(\varphi_i)d\mu_i
\end{align*}	
	\end{lemma}
	
	\begin{proof} 
Since $F_i$ is supperlinear one has $\mathring{D}(F_i^\circ)=  \mathbb{R} $ for $i=1,2$. Let any $(\varphi_1,\varphi_2)\in \Lambda$ then $(\varphi_1,\varphi_2)\in C_b(X_1)\times C_b(X_2)$ and $\varphi_1(x_1)\leq C(x_1,p)+p(-\varphi_2)$ for every $x_1\in X_1, p\in\cP(X_2)$. This implies that $\varphi_1(x_1)\leq R_C(-\varphi_2)(x_1)$ for every $x_1\in X_1$. Moreover, from \eqref{temp 10} one gets $F_i^\circ$ is also nondecreasing on $(-(F_i)'_\infty,+\infty)=  \mathbb{R} $ for $i=1,2$. Therefore, \begin{align*}
\sum_{i=1}^2\int_{X_i}F_i^\circ(\varphi_i)d\mu_i&\leq \int_{X_1}F_1^\circ(R_C(-\varphi_2))d\mu_1+\int_{X_2}F_2^\circ(\varphi_2)d\mu_2\\
&\leq \sup_{\varphi\in C_b(X_2)}\int F_1^\circ (R_C \varphi)d\mu_1+\int F_2^\circ (-\varphi) d\mu_2. 
\end{align*} 		
So that we obtain \begin{align*}
\sup_{(\varphi_1,\varphi_2)\in \Lambda}\sum_{i=1}^2\int_{X_i}F_i^\circ(\varphi_i)d\mu_i\le\sup_{\varphi\in C_b(X_2)}\int F_1^\circ (R_C \varphi)d\mu_1+\int F_2^\circ (-\varphi) d\mu_2.
\end{align*}
Now we prove that $$\cE_C(\mu_1,\mu_2)\geq \sup_{\varphi\in C_b(X_2)}\int_{X_1}F_1^\circ(R_C\varphi)d\mu_1+\int_{X_2}F_2^\circ(-\varphi)d\mu_2.$$
If the problem (WOET) is not feasible then it is clear as $\cE_C(\mu_1,\mu_2)=+\infty$. Now we assume the feasibility of the problem (WOET). Applying Theorem \ref{T-existence of minimizers}, there exist minimizers of the problem (WOET). Let $ \boldsymbol\gamma\in \cM(X_1\times X_2) $ be an optimal plan for problem (WOET). We will show that  $ R_C\varphi\in L^1(X_1,\gamma_1) $ for every $\varphi\in C_b(X_2)$. Since $ \varphi\in C_b(X_2) $ and $C$ is bounded from below we can assume there are $ M_1,M_2>0$ such that $\varphi>-M_1$ and $C(x_1,p)> -M_2$ for any $x_1\in X_1,p\in\cP(X_2)$. Then one has $ R_C\varphi\geq -M_1-M_2$ on $X_1$. Thus, $  \vert R_C\varphi(x_1) \vert \le \max\{M_1+M_2,C(x_1,\gamma_{x_1})+\gamma_{x_1}(\varphi) \}, $ for every $x_1\in X_1$. On the other hand, we have
\begin{align*}
 \int_{X_1}[C(x_1,\gamma_{x_1})+\gamma_{x_1}(\varphi)]d\gamma_1(x_1)\leq \cE_C(\boldsymbol{\gamma} \vert \mu_1,\mu_2)+M \vert \gamma_1 \vert <\infty. 
\end{align*} 
 Hence $ R_C\varphi\in L^1(X_1,\gamma_1).$

Let any $\varphi\in C_b(X_2)$, as $R_C\varphi(x_1)+p(-\varphi)\leq C(x_1,p)$ for every $x_1\in X_1,p\in \cP(X_2)$ one has 
\begin{eqnarray*}
			\cE_C(\boldsymbol{\gamma} \vert \mu_1,\mu_2)&=&\cF_1(\gamma_1 \vert \mu_1)+\cF_2(\gamma_2 \vert \mu_2)+\int_{X_1}C(x_1,\gamma_{x_1})d\gamma_1(x_1)\\
			&\geq & \cF_1(\gamma_1 \vert \mu_1)+\cF_2(\gamma_2 \vert \mu_2)+\int_{X_1}(R_C\varphi(x_1)+\gamma_{x_1}(-\varphi))d\gamma_1(x_1)\\
			&=& \cF_1(\gamma_1 \vert \mu_1)+\cF_2(\gamma_2 \vert \mu_2)+\int_{X_1}R_C\varphi d\gamma_1+\int_{X_1}\int_{X_2}(-\varphi)(x_2)d\gamma_{x_1}d\gamma_1\\
			&=& \cF_1(\gamma_1 \vert \mu_1)+\int_{X_1}R_C\varphi d\gamma_1+\cF_2(\gamma_2 \vert \mu_2)+\int_{X_2}(-\varphi)d\gamma_2.
		\end{eqnarray*} 
Since $\cF(\gamma_1 \vert \mu_1)<\infty$ and $R_C\varphi\in L^1(\gamma_1)$ applying Lemma \ref{L-equality of F for (BM)} ($\psi=F_1^\circ(R_C\varphi), \phi=-R_C\varphi$) we obtain that $$\cF_1(\gamma_1 \vert \mu_1)+\int_{X_1}R_C\varphi d\gamma_1\geq \int_{X_1}F_1^\circ(R_C\varphi)d\mu_1.$$
Similarly, we also have that $$\cF_2(\gamma_2 \vert \mu_2)+\int_{X_2}(-\varphi)d\gamma_2\geq \int_{X_2}F_2^\circ(-\varphi)d\mu_2.$$
Therefore, $$\cE_C(\mu_1,\mu_2)\geq \int_{X_1}F_1^\circ(R_C\varphi)d\mu_1+\int_{X_2}F_2^\circ(-\varphi)d\mu_2.$$
Applying Lemma \ref{L-dual formula for Lamda set} we get that 
\begin{align*}
\cE_C(\mu_1,\mu_2)&=\sup_{(\varphi_1,\varphi_2)\in \Lambda}\sum_{i=1}^2\int_{X_i}F_i^\circ(\varphi_i)d\mu_i\\
&=\int_{X_1}F_1^\circ(R_C\varphi)d\mu_1+\int_{X_2}F_2^\circ(-\varphi)d\mu_2.
\end{align*} 
	   \end{proof}
\begin{proof}[Proof of Theorem \ref{dual theorem}] Theorem \ref{dual theorem} follows from Lemmas \ref{L-dual theorem} and \ref{reduce to RC}.
   \end{proof}

	Next, we want to investigate the monotonicity property of the optimal plans of problem (WOET). 
	\begin{definition}
	(\cite[Definition 5.1]{BBP})
		We say that a measure $ \boldsymbol{\gamma}\in\cM(X_1\times X_2) $ is $ C $-monotone if there exists a measurable set $ \Gamma\subseteq X_1 $ such that $ \gamma_1 $ is concentrated on $ \Gamma $ and for any finite number of points $ x_1^1,\dots,x_1^N $ in $ \Gamma $, for any measures $ m_1,\dots,m_N $ in $ \cP(X_2) $ with $ \sum_{i=1}^Nm_i=\sum_{i=1}^N\gamma_{x_1^i} $, the follow inequality holds:\[ \sum_{i=1}^NC(x_1^i,\gamma_{x_1^i})\le\sum_{i=1}^NC(x_1^i,m_i). \]
	\end{definition}
	\begin{corollary}
		Assume that problem (WOET) is feasible and coercive for $\mu_i\in \cM(X_i),i=1,2$. If $\boldsymbol{\gamma}\in \cM(X_1\times X_2) $ is an optimal plan for $\cE_C(\mu_1,\mu_2)$ then $\boldsymbol{\gamma}$ is $C$-monotone.
	\end{corollary}
	\begin{proof}
		The case that $ \boldsymbol{\gamma} $ is the null measure is a trivial case so we can assume $ \boldsymbol{\gamma} $ is not the null measure. Because $ \boldsymbol{\gamma} $ is an optimal plan for the problem (WOET) we get that $ \boldsymbol{\gamma}/ \vert \boldsymbol{\gamma} \vert  \in \cP(X_1\times X_2)$ is an optimal plan for weak transport costs problem for its marginals discussed in \cite{BBP}. Applying \cite[Theorem 5.3]{BBP} we get the result.
	   \end{proof}
	
		\section{Examples}
	In this section, we will illustrate examples of entropy functions $F_i$ and cost functions $C:X_1\times \cP(X_2)\to (-\infty,+\infty]$ for our Weak Optimal Entropy Transport Problems. 
	\begin{example}
	\label{E-Optimal Entropy-Transport Problems}
	(Optimal Entropy-Transport Problems)
	If there exists some cost function $c:X_1\times X_2\to (-\infty,+\infty]$ which is lower semi-continuous and bounded from below, such that $C(x_1,p)=\int_{X_2}c(x_1,x_2)dp(x_2)$ for every $x_1\in X_1,p\in \cP(X_2)$ then the problem (WOET) becomes the Optimal-Entropy Transport problem \cite{LMS18} of finding $\bar{\boldsymbol{\gamma}}\in \cM(X_1\times X_2)$ minimizing
	$$\cE(\bar{\boldsymbol{\gamma}} \vert \mu_1,\mu_2)=\cE(\mu_1,\mu_2):=\inf_{\boldsymbol{\gamma}\in \cM(X_1\times X_2)} \cE(\boldsymbol{\gamma} \vert \mu_1,\mu_2), $$
		where $\cE(\boldsymbol{\gamma} \vert \mu_1,\mu_2):=\sum_{i=1}^2\cF_i(\gamma_i \vert \mu_i)+\int_{X_1\times X_2}c(x_1,x_2)d\boldsymbol{\gamma}(x_1,x_2)$.

Since $c$ is bounded from below, so is $C$. Moreover, applying the following lemma we get that $C$ is lower semi-continuous on $X_1\times\cP(X_2)$.
\begin{lemma}\label{L-lower semicontinuous in Martingale problems}
Let $X_1,X_2$ be Polish metric spaces and let $f:X_1\times X_2\to (-\infty,+\infty]$ be a lower semi-continuous function satisfying that $f$ is bounded from below. Let $(x^n,p^n)\subset X_1\times\cP(X_2)$ such that $(x^n,p^n)\to (x^0,p^0)$ as $n\to \infty$, for $(x^0,p^0)\in X_1\times \cP(X_2)$. Then we have $$\liminf_{n\to\infty}\int_{X_2}f(x^n,x_2)dp^n(x_2)\geq \int_{X_2}f(x^0,x_2)dp^0(x_2).$$
\end{lemma}
\begin{proof}
For any $n\in  \mathbb{N} $, we define $\boldsymbol{P}^n:=\delta_{x^n}\otimes p^n\in\cP(X_1\times X_2)$ and set $\boldsymbol{P}^0:=\delta_{x^0}\otimes p^0\in\cP(X_1\times X_2)$. Since $\lim_{n\to\infty}x^n=x^0$ one gets that $\delta_{x^n}\to \delta_{x^0}$ as $n\to \infty$ under the weak topology. Hence, by \cite[Theorem 2.8 (ii)]{B} we obtain that $\boldsymbol{P}^n\to \boldsymbol{P}^0$ as $n\to \infty$ under the weak topology. Moreover, as $f$ is lower semi-continuous and bounded from below we get that \begin{align*}
\liminf_{n\to\infty}\int_{X_2}f(x^n,x_2)dp^n(x_2)&=\liminf_{n\to\infty}\int_{X_1\times X_2}f(x_1,x_2)d\boldsymbol{P}^n(x_1,x_2)\\
&\geq \int_{X_1\times X_2}f(x_1,x_2)d\boldsymbol{P}^0(x_1,x_2)\\
&=\int_{X_2}f(x^0,x_2)dp^0(x_2).
\end{align*}
Hence, we get the result.
   \end{proof}
   \begin{lemma}
   \label{L-temp1}
   Let $X_1$ and $X_2$ be Polish metric spaces. Let $F_1, F_2: [0,\infty)\to [0,\infty]$ be admissible entropy functions. Assume that there exists some function $c:X_1\times X_2\to (-\infty,+\infty]$ which is lower semi-continuous and bounded from below, such that $C(x_1,p)=\int_{X_2}c(x_1,x_2)dp(x_2)$ for every $x_1\in X_1,p\in \cP(X_2)$. We recall
   \begin{align*}
	\Lambda:=\bigg\{(\varphi_1,\varphi_2)\in C_b(X_1,\mathring{ D}(F_1^\circ))\times C_b(X_2,\mathring{ D}(F_2^\circ))&:\varphi_1(x_1)+p(\varphi_2)\leq C(x_1,p),\\
 &\mbox{ for every }
	 x_1\in X_1, p\in \cP(X_2)\bigg\}.
\end{align*}
\begin{align*}
\boldsymbol{\Phi}&:=\bigg\{(\varphi_1,\varphi_2)\in C_b(X_1,\mathring{D}(F_1^\circ))\times C_b(X_2,\mathring{D}(F_2^\circ)):\varphi_1\oplus \varphi_2\leq c\bigg\}.
\end{align*}
Then $\Lambda=\boldsymbol{\Phi}$. 
   \end{lemma}
\begin{proof}

Let $(\varphi_1,\varphi_2)\in \Lambda$. Then for every $x_1\in X_1, x_2\in X_2$ we have that $$\varphi_1(x_1)+\varphi_2(x_2)=\varphi_1(x_1)+\delta_{x_2}(\varphi_2)\leq \int_{X_2}c(x_1,x'_2)\delta_{x_2}(x_2')=c(x_1,x_2).$$ Therefore 
$\Lambda\subset \boldsymbol{\Phi}$. 

Conversely, let $(\varphi_1,\varphi_2)\in \boldsymbol{\Phi}$. For every $x_1\in X_1$ and $p\in \cP(X_2)$ we have that $\varphi_1(x_1)+p(\varphi_2)=\int_{X_2}(\varphi_1(x_1)+\varphi_2(x_2))dp(x_2)\leq \int_{X_2}c(x_1,x_2)dp(x_2)=C(x_1,p).$ Hence $\boldsymbol{\Phi}\subset \Lambda$. 
\end{proof}

	\end{example}
	\begin{example}
	\label{E-Weak Optimal Transport Problems}
	(Weak Optimal Transport Problems)
	For $i\in\{1,2\}$, we define the admissible entropy functions $F_i:[0,\infty)\rightarrow [0,\infty]$ by $$F_i(r):=\left\lbrace\begin{array}{ll}
		0 & \text{ if } r=1,\\
		+\infty &\text{ otherwise.}
	\end{array}\right.$$	
	Then the problem $(WOET)$ becomes the pure weak transport problem 
	\begin{equation}
\cE_C(\mu_1,\mu_2)=\inf_{\boldsymbol{\gamma}\in \cM(X_1\times X_2)}\left\{\int_{X_1}C(x_1,\gamma_{x_1})d\gamma_1(x_1) \vert \pi^i_\sharp\gamma=\mu_i,i=1,2\right\}.\label{P-temp}
	\end{equation}
		In this example, if $\boldsymbol{\gamma}\in \cM(X_1\times X_2)$ is a feasible plan then $\mu_1,\mu_2$ are the marginals of $\boldsymbol{\gamma}$. Thus, a necessary condition for feasibility is that $\vert \mu_1\vert=\vert\mu_2\vert$. If furthermore $\mu_i\in \cP(X_i), i=1,2$ then \eqref{P-temp} will be the weak transport problem which has been introduced by \cite{GRST}. 

	In addition to, if $X_1=X_2=X\subset  \mathbb{R} ^d$ for some $d\in   \mathbb{N} $ and
			\[  C(x_1,p)=\begin{cases} 
				\int_{X}c(x_1,x_2)dp(x_2) & \mbox{ if } \int_X x_2dp(x_2)=x_1, \\
				+\infty & \mbox{ otherwise },
				
			\end{cases}
			\] 
			then the problem \eqref{P-temp} will become the classical martingale optimal transport problem for every $\mu_1,\mu_2\in \cP(X)$. It was introduced first for the case $X=\mathbb{R}$ by Beiglb\"{o}ck, Henry-Labord\`ere and Penkner \cite{BHP} and since then it has been studied intensively \cite{BJ,BP,BBHK,GKL,HT}. Now we introduce our martingale optimal entropy transport (MOET) problems. Given $\mu,\nu\in \cM(X)$, we denote by $\Pi_M(\mu,\nu)$ the set of all measures $\boldsymbol{\gamma}\in \cM(X^2)$ such that $\pi^1_\sharp \boldsymbol{\gamma}=\mu, \pi^2_\sharp \boldsymbol{\gamma}=\nu$ and $\int_X yd\pi_x(y)=x \mbox{ } \mu$-almost everywhere, where $(\pi_x)_{x\in X}$ is the disintegration of $\boldsymbol{\gamma}$ with respect to $\mu$. We denote by $\cM_M(X^2)$ the set of all $\boldsymbol{\gamma}\in \cM(X^2)$ such that $\boldsymbol{\gamma}\in \Pi_M(\pi^1\gamma,\pi^2\gamma)$. Our (MOET) problem is defined as
\begin{align*}
			\cE_M(\mu_1,\mu_2):=\inf_{\boldsymbol{\gamma}\in \cM(X^2)}\cE_C(\boldsymbol{\gamma} \vert \mu_1,\mu_2)=\inf_{\boldsymbol{\gamma}\in \cM_M(X^2)}\left\{\sum_{i=1}^2\cF(\gamma_i \vert \mu_i)+\int_{X\times X}c(x_1,x_2)d\boldsymbol{\gamma}\right\}.
		\end{align*}	
	\end{example}	
Using the ideas of \cite[Section 5]{LMS18}, we can establish a Kantorovich duality of our (MOET) problem in terms of homogeneous marginal perspective functionals and homogeneous constraints. However, as we have not found its applications yet, we skip the details here.
	\begin{example}
	(Weak Logarithmic Entropy Transport (WLET))	
	Suppose that $X_1=X_2=X$ is a Polish space and let $F_i(t)=t\log t-t+1$ for $t\geq 0$, $i=1,2$ with the convention that $0\log 0=0$. This entropy functional plays an important role in the study of Optimal Entropy Transport problems \cite[Sections 6-8]{LMS18}.
	In this case, $F_i$ is superlinear and hence our (WOET) problem becomes the Weak Logarithmic Entropy Transport problem \begin{align*}
	\cE(\mu_1,\mu_2)&=\WLET(\mu_1,\mu_2)\\
	&=\inf_{\boldsymbol{\gamma}\in \cM(X\times X)}\left\{\sum_{i=1}^2\int_X(\sigma_i\log\sigma_i-\sigma_i+1)d\mu_i+\int_XC(x_1,\gamma_{x_1})d\gamma_1(x_1)\right\},
\end{align*}	 
where $\sigma_i=\dfrac{d\gamma_i}{d\mu_i}$.
	
	The feasible condition always holds from Lemma \ref{L-feasible} since $F_1(0)=F_2(0)=1<\infty$. Furthermore, $R_i(r)=rF_i(1/r)=r-1-\log r$ for $r>0$ and $R_i(0)=+\infty$; and $R^*_i(\psi)=-\log(1-\psi)$ for $\psi<1$ and $R^*_i(\psi)=+\infty$ for $\psi\geq 1$. 

	\end{example}

\begin{example}(The $\chi^2$-divergence)
In this example, let $F_1\in \Adm(  \mathbb{R} _+)$ such that $F_1$ is superlinear and $F_1(0)<\infty$. We consider $F_2(t)=\varphi_{\chi^2}(t)=(t-1)^2$ for $t\geq 0$. As $F_1(0)<\infty$ and $F_2(0)=1$ one has that the problem (EWOT) is feasible. Observe that $F_2$ is superlinear and $R_2(r)=(r-1)^2/r$ for $r>0$ and $R_2(0)=+\infty$. From this, it is not difficult to check that $$R^*_2(\psi)=\sup_{r\geq 0}\{r\psi-R_2(r)\}=\left\lbrace\begin{array}{ll}
		+\infty &\text{ if } \psi>1,\\
		2-2\sqrt{1-\psi} &\text{ if } \psi\leq 1.
	\end{array}\right.$$
\end{example}
\begin{example}
	(Marton's cost functions)
	Let $X$ be a compact subset of $  \mathbb{R} ^m$ and let $C:X\times \cP(X)\to (-\infty,+\infty]$ be the cost function defined by $$C(x,p):=\theta\bigg(x-\int_X ydp(y)\bigg), \mbox{ for every } x\in X, p\in \cP(X),$$
where $\theta:  \mathbb{R} ^m\to (-\infty,+\infty]$ is a lower semi-continuous convex function such that $\theta$ is bounded from below. Then $C(x,\cdot)$ is convex on $\cP(X)$ for every $x\in X$ and $C$ is bounded from below. Next, we will check that $C$ is lower semi-continuous on $X\times\cP(X)$. Let $\{(x_n,p_n)\}_n\subset X\times \cP(X)$ such that $(x_n,p_n)\to (x_0,p_0)$ as $n\to \infty$ for $(x_0,p_0)\in X\times \cP(X)$. As $X$ is compact, one gets that $$\lim_{n\to\infty}\left(x_n-\int_Xydp_n(y)\right)=x_0-\int_Xydp_0(y).$$
Moreover, since $\theta$ is lower semi-continuous we obtain that $$\liminf_{n\to\infty}C(x_n,p_n)=\liminf_{n\to\infty}\theta\left(x_n-\int_Xydp_n(y)\right)\geq \theta\left(x_0-\int_Xydp_0(y)\right)=C(x_0,p_0).$$
This means that $C$ is lower semi-continuous on $X\times\cP(X)$.

The following theorem is an extension of \cite[Theorem 2.11]{GRST}.		
\begin{theorem}
Let $X$ be a compact, convex subset of $  \mathbb{R} ^m$. Assume that $F_1,F_2$ are superlinear. For every  $\mu_1,\mu_2\in \cM(X)$ we have
\begin{align*}
			\cE_C(\mu_1,\mu_2)=\sup\bigg\{ \int_X F_1^\circ (R_\theta \varphi)d\mu_1+\int_X F_2^\circ (-\varphi) d\mu_2: \varphi\in LSC_{bc}(X) \bigg\}.
		\end{align*}	
where $R_\theta\varphi (x):=\inf_{p\in \cP(X)}\{C(x,p)+p(\varphi)\}$ and $LSC_{bc}(X)$ is the set of all bounded, lower semi-continuous and convex function on $X$.
	\end{theorem}
	\begin{proof}
For every $p\in \cP(X)$ we will show that $\int_{X}ydp(y)\in X$. If $p=\sum_{i=1}^N \lambda_i\delta_{x_i}$ where $\sum_{i=1}^N\lambda_i=1$ and $x_i\in X$ for $i=1,\ldots,N$ then as $X$ is convex one has $$\int_Xydp(y)=\sum_{i=1}^N\lambda_ix_i\in X.$$
Now, let any $p\in\cP(X)$. As $X$ is compact, applying \cite[Theorem 5.9]{San} and \cite[Theorem 6.18]{V}, we can approximate $p$ by a sequence of probability measures with finite support in the weak topology. Thus, since $X$ is closed we get that $\int_Xydp(y)\in X$.

For any $\varphi\in C_b(X)$, we define the function $g_\varphi:X\to   \mathbb{R} $	by $$g_\varphi(z):=\inf_{p\in \cP(X)}\{\int_{X}\varphi dp:\int_{X} yp(dy)=z \},$$
for every $z\in X$. Then it is not difficult to check that $g$ is convex on $X$. Since $\varphi\in C_b(X)$, there exists $m\in   \mathbb{R} $ such that $\varphi(x)\geq m$ for every $x\in X$. Then for any $p\in \cP(X)$ we have  $\int_{X}\varphi(y)dp(y)\geq m.$ So $g_\varphi(z)\geq m$ for every $z\in X$. Furthermore, for every $z\in X$ one has	
$$g_\varphi(z)\leq \int_{X}\varphi d\delta_z=\varphi(z).$$ So $g_\varphi$ is bounded. Next, we will check that $g_\varphi$ is the greatest convex function bounded above by $\varphi$. Let any convex function $\widehat{\varphi}$ such that $m\leq \widehat{\varphi}(x)\leq \varphi(x)$ for every  $x\in X$. Then for any $z\in X$ let $p\in\cP(X)$ such that $\int_Xydp(y)=z$, applying Jensen's inequality for $\widehat{\varphi}$  
one has $$\int_{X}\varphi(y) dp(y)\geq \int_X\widehat{\varphi}(y)dp(y)\geq \widehat{\varphi}(\int_Xydp(y))=\widehat{\varphi}(z).$$	
Hence, $g_\varphi\geq \widehat{\varphi}$ on $X$. This means that $g_\varphi$ is the greatest convex function bounded above by $\varphi$. Now, we extend the function $\varphi$ by putting $\varphi(x)=+\infty$ for every $x\notin X$. 
Then by \cite[Corollary 17.2.1]{R} we obtain that $g_\varphi$ is lower semi-continuous on $X$.

On the other hand, by the definition of $g_\varphi$, for every $x\in X$, we get that \begin{align*}
R_{\theta}\varphi(x)&=\inf_{p\in\cP(X)}\bigg\{\int_X\varphi dp+\theta\bigg(x-\int_Xydp\bigg)\bigg\}\\
&=\inf_{z\in X}\bigg\{g_\varphi(z)+\theta(x-z)\bigg\}.
\end{align*}		
Furthermore, we have $\inf_{z\in X}\{g_\varphi(z)+\theta(x-z)\}\leq R_\theta g_\varphi(x)$ for every $x\in X$. Indeed, for any $p\in \cP(X)$ setting $w:=\int_Xydp(y)\in X$. For every $x\in X$, using Jensen's inequality again for the convex function $g_\varphi$ we get \begin{align*}
\int_Xg_\varphi dp+\theta\left(x-\int_Xydp\right)\geq g_\varphi(w)+\theta(x-w)\geq \inf_{z\in X}\{g_\varphi(z)+\theta(x-z)\}.
\end{align*}
Combining with $g_\varphi\leq \varphi$ on $X$, one gets that 		
$$R_\theta\varphi(x)=\inf_{z\in X}\{g_\varphi(z)+\theta(x-z)\}\leq R_\theta g_\varphi(x)\leq R_\theta\varphi (x).$$		
		 Hence from \eqref{temp 10} we get
		  \begin{align*} \int_X F_1^\circ (R_\theta \varphi)d\mu_1+\int_X F_2^\circ (-\varphi) d\mu_2=&\int_X F_1^\circ (R_\theta g_\varphi)d\mu_1+\int_X F_2^\circ (-\varphi) d\mu_2\\
		  \le&\int_X F_1^\circ (R_\theta g_\varphi)d\mu_1+\int_X F_2^\circ (-g_\varphi) d\mu_2. \end{align*}
		Therefore, applying Lemma \ref{reduce to RC} we obtain
		\begin{align*}
			\cE_C(\mu_1,\mu_2)=&\sup\bigg\{ \int_X F_1^\circ (R_\theta \varphi)d\mu_1+\int_X F_2^\circ (-\varphi) d\mu_2: \varphi\in C_b(X)\bigg\}\\&\le\sup\bigg\{ \int_X F_1^\circ (R_\theta \varphi)d\mu_1+\int_X F_2^\circ (-\varphi) d\mu_2: \varphi\in LSC_{bc}(X)\bigg\}.
		\end{align*}
To complete the proof, we only need to prove that \begin{align}\label{F-marton cost}
\cE_C(\mu_1,\mu_2)\geq \sup\bigg\{ \int_X F_1^\circ (R_\theta \varphi)d\mu_1+\int_X F_2^\circ (-\varphi) d\mu_2: \varphi\in LSC_{bc}(X) \bigg\}.
\end{align}	
If the problem (WOET) is not feasible then both sides of \eqref{F-marton cost} are infinity. So we can assume the problem (WOET) is feasible. By Theorem \ref{T-existence of minimizers}, let $\boldsymbol{\gamma}\in \cM(X\times X)$ such that $\cE_C(\mu_1,\mu_2)=\cE_C(\boldsymbol{\gamma} \vert \mu_1,\mu_2)$. For every $\varphi\in LSC_{bc}(X)$, we have \begin{align*}
\cE_C(\boldsymbol{\gamma} \vert \mu_1,\mu_2)=&\cF_1(\gamma_1 \vert \mu_1)+\cF_2(\gamma_2 \vert \mu_2)+\int_{X}C(x_1,\gamma_{x_1})d\gamma_1(x_1)\\
			\geq & \cF_1(\gamma_1 \vert \mu_1)+\cF_2(\gamma_2 \vert \mu_2)+\int_{X}(R_\theta\varphi(x_1)+\gamma_{x_1}(-\varphi))d\gamma_1(x_1)\\
			=& \cF_1(\gamma_1 \vert \mu_1)+\int_{X}R_\theta\varphi d\gamma_1+\cF_2(\gamma_2 \vert \mu_2)+\int_{X}(-\varphi)d\gamma_2.
\end{align*}
Since $\varphi$ is bounded, using the same arguments as in the proof of Lemma \ref{reduce to RC} one has $R_\theta\varphi\in L^1(X,\gamma_1)$. Hence, by Lemma \ref{L-equality of F for (BM)} one gets 
$$\cF_1(\gamma_1 \vert \mu_1)+\int_{X}R_\theta\varphi d\gamma_1\geq \int_{X}F_1^\circ(R_\theta\varphi)d\mu_1,$$
$$\cF_2(\gamma_2 \vert \mu_2)+\int_{X}(-\varphi)d\gamma_2\geq \int_{X_2}F_2^\circ(-\varphi)d\mu_2.$$
This implies that \eqref{F-marton cost} and then we get the result.
	   \end{proof}
	
	\end{example}
	\section{Declarations}
	~~~~~\textbf{Ethical Approval} (applicable for both human and/ or animal studies. Ethical committees, Internal Review Boards and guidelines followed must be named. When applicable, additional headings with statements on \textbf{consent to participate} and \textbf{consent to publish} are also required): not applicable.
	 
	   \textbf{Competing interests} (always applicable and includes interests of a financial or personal nature) 
	 
	 The authors declare no competing interests.	   
	
   \textbf{Authors' contributions} (applicable for submissions with multiple authors) 
    
    Authors contributed equally.
    
     \textbf{Funding} (details of any funding received) 
	 
	  This work was supported by the National Research Foundation of Korea (NRF) grant funded by the Korea government (MSIT) No. NRF-2016R1A5A1008055 and No. NRF-2019R1C1C1007107. 		
    
     \textbf{Availability of data and materials} (a statement on how any datasets used can be accessed): not applicable.
	


\begin{thebibliography}{9}
\bibitem{ABC} J.-J. Alibert, G. Bouchitt\'{e}, T. Champion, A new class of costs for optimal transport planning, European J. Appl. Math., 30(6) (2019) 1229-1263.

\bibitem{ACJ}  Aur\'{e}lien Alfonsi, Jacopo Corbetta, Benjamin Jourdain, Sampling of probability measures in the convex order by Wasserstein projection, English, with English and French summaries, Ann. Inst. Henri Poincar\'{e} Probab. Stat., 56(3) (2020) 1706-1729.

\bibitem{AGS} Luigi Ambrosio, Nicola Gigli, Giuseppe Savar\'{e}, Gradient flows in metric spaces and in the space of probability measures, Lectures in Mathematics ETH Z\"{u}rich, Birkh\"{a}user Verlag, Basel (2005).

\bibitem{BBHK} Julio Backhoff-Veraguas, Mathias Beiglb\"{o}ck, Martin Huesmann, Sigrid K\"{a}llblad, Martingale Benamou-Brenier: a probabilistic perspective, Ann. Probab. 48 (5) (2020), 2258–2289.

\bibitem{BBP} J. Backhoff-Veraguas, M. Beiglb\"{o}ck, and G. Pammer, Existence, duality, and cyclical monotonicity for weak transport costs, Calc. Var. Partial Differential Equations 58(6) (2019), Paper No. 203, 28.

\bibitem{BP} J. Backhoff-Veraguas, G. Pammer, Stability of martingale optimal transport and weak optimal transport, Ann. Appl. Probab. 32 (1) (2022), 721–752.

\bibitem{BHP} Mathias Beiglb\"{o}ck, Pierre Henry-Labord\`ere, Friedrich Penkner, Model-independent bounds for option prices—a mass transport approach, Finance Stoch. 17 (3) (2013), 477–501.

\bibitem{BJ} Mathias Beiglb\"{o}ck, Nicolas Juillet, On a problem of optimal transport under marginal martingale constraints, Ann. Probab. 44  (1) (2016), 42–106.

\bibitem{B} Patrick Billingsley, Convergence of probability measures, 2nd ed., Wiley Series in Prob- ability and Statistics: Probability and Statistics, John Wiley \& Sons, Inc., New York, 1999. A Wiley-Interscience Publication.

\bibitem{CP} Nhan-Phu Chung, Minh-Nhat Phung, Barycenters in the Hellinger-Kantorovich Space, Appl. Math. Optim. 84 (2) (2021), 1791–1820.

\bibitem{CT} Nhan-Phu Chung, Thanh-Son Trinh, Duality and quotient spaces of generalized Wasserstein spaces, arXiv:1904.12461.

\bibitem{CT1} Nhan-Phu Chung, Thanh-Son Trinh, Unbalanced optimal total variation transport problems and generalized Wasserstein barycenters, Proc. Roy. Soc. Edinburgh Sect. A 152 (3) (2022), 674–700.

\bibitem{D} Nicol\`o De Ponti, Metric properties of homogeneous and spatially inhomogeneous F-
divergences, IEEE Trans. Inform. Theory 66 (5) (2020), 2872–2890.

\bibitem{ET} Ivar Ekeland and Roger T\'{e}mam, Convex analysis and variational problems, Corrected reprint of the 1976 English edition, Classics in Applied Mathematics, vol. 28, Society for Industrial and Applied Mathematics (SIAM), Philadelphia, PA, 1999. Translated
from the French.

\bibitem{FMS}Gero Friesecke, Daniel Matthes, Bernhard Schmitzer, Barycenters for the Hellinger-Kantorovich distance over $\mathbb{R}^d$, SIAM J. Math. Anal. 53 (1) (2021), 62–110.

\bibitem{Galichon} Alfred Galichon, Optimal transport methods in economics, Princeton University Press,
Princeton, NJ, 2016.

\bibitem{GKL} Nassif Ghoussoub, Young-Heon Kim, and Tongseok Lim, Structure of optimal martingale transport plans in general dimensions, Ann. Probab. 47 (1) (2019), 109–164.

\bibitem{GT} Nicola Gigli and Luca Tamanini, Benamou-Brenier and duality formulas for the entropic cost on ${RCD}^*(K,N)$ spaces, Probab. Theory Related Fields 176 (1-2) (2020), 1–34.

\bibitem{GJ} Nathael Gozlan, Nicolas Juillet, On a mixture of Brenier and Strassen theorems,
Proc. Lond. Math. Soc. 120 (3) (2020), 434–463.

\bibitem{GRST} Nathael Gozlan, Cyril Roberto, Paul-Marie Samson, Prasad Tetali, Kantorovich
duality for general transport costs and applications, J. Funct. Anal. 273 (11) (2017),
3327–3405.

\bibitem{GRSST} Nathael Gozlan, Cyril Roberto, Paul-Marie Samson, Yan Shu, Prasad Tetali, Characterization of a class of weak transport-entropy inequalities on the line, Ann. Inst. Henri Poincar\'{e} Probab. Stat. 54 (3) (2018), 1667–1693 (English, with English and French summaries).

\bibitem{HT} Martin Huesmann, Dario Trevisan, A Benamou-Brenier formulation of martingale optimal transport, Bernoulli 25 (4A) (2019), 2729–2757.

\bibitem{Kant42} L. Kantorovitch, On the translocation of masses, C. R. (Doklady) Acad. Sci. URSS (N.S.) 37 (1942), 199–201.

\bibitem{Kant48} L. V. Kantorovich, On a problem of Monge, Zap. Nauchn. Sem. S.-Peterburg. Otdel. Mat. Inst. Steklov. (POMI) 312 (2004), no. Teor. Predst. Din. Sist. Komb. i Algoritm. Metody. 11, 15–16 (Russian); English transl., J. Math. Sci. (N.Y.) 133 (4) (2006), 1383.

\bibitem{KV} Stanislav Kondratyev, Dmitry Vorotnikov, Convex Sobolev inequalities related to unbalanced optimal transport, J. Differential Equations 268 (7) (2020), 3705–3724.

\bibitem{LM} Vaios Laschos, Alexander Mielke, Geometric properties of cones with applications on the Hellinger-Kantorovich space, and a new distance on the space of probability measures, J. Funct. Anal. 276 (11) (2019), 3529–3576.

\bibitem{LMS18} Matthias Liero, Alexander Mielke, Giuseppe Savar\'{e}, Optimal entropy-transport problems and a new Hellinger-Kantorovich distance between positive measures, Invent. Math. 211 (3) (2018), 969–1117.

\bibitem{Par} K. R. Parthasarathy, {\it Probability measures on metric spaces}, AMS Chelsea Publishing, Providence, RI, 2005. Reprint of the 1967 original.

\bibitem{PC} Gabriel Peyr\'{e}, Marco Cuturi, Computational Optimal Transport, Foundations and Trends in Machine Learning 11 (5-6) (2019), 355–607.

\bibitem{MS} Matteo Muratori, Giuseppe Savar\'{e}, Gradient flows and evolution variational in- equalities in metric spaces. I: Structural properties, J. Funct. Anal. 278 (4) (2020), 108347, 67.

\bibitem{R} R. Tyrrell Rockafellar, Convex analysis, Princeton Mathematical Series, No. 28, Prince- ton University Press, Princeton, N.J., 1970.

\bibitem{San} Filippo Santambrogio, Optimal transport for applied mathematicians, Progress in Non- linear Differential Equations and their Applications, vol. 87, Birkh\"{a}user/Springer, Cham, 2015. Calculus of variations, PDEs, and modeling.

\bibitem{Shu} Yan Shu, From Hopf-Lax formula to optimal weak transfer plan, SIAM J. Math. Anal. 52 (3) (2020), 3052–3072.

\bibitem{V03} C\'{e}dric Villani, Topics in optimal transportation, Graduate Studies in Mathematics, vol. 58, American Mathematical Society, Providence, RI, 2003.

\bibitem{V} C\'{e}dric Villani, Optimal transport, Grundlehren der Mathematischen Wissenschaften [Funda- mental Principles of Mathematical Sciences], vol. 338, Springer-Verlag, Berlin, 2009. Old and new.


\end{thebibliography}

%


\end{document}